\documentclass[12pt]{amsart}
\usepackage{latexsym, amssymb, amscd, rotating}
\renewcommand\a{\alpha}
\renewcommand\b{\beta}
\newcommand\g{\gamma}
\renewcommand\d{\delta}
\newcommand\la{\lambda}
\newcommand\z{\zeta}
\newcommand\e{\eta}
\renewcommand\th{\theta}

\newcommand\s{\sigma}
\newcommand\x{\chi}
\newcommand\f{\phi}
\newcommand\vf{\varphi}
\newcommand\p{\psi}
\renewcommand\t{\tau}
\renewcommand\r{\rho}

\newcommand\Om{\Omega}
\newcommand\w{\omega}

\newcommand\vD{\varDelta}

\newcommand\vL{\varLambda}

\newcommand\vG{\varGamma}
\newcommand\ve{\varepsilon}

\newcommand{\FF}{\mathbb F}

\newcommand\Fq{{\mathbf F}_q}
\newcommand\Fqr{{\mathbf F}_{q^r}}
\newcommand\Fqq{{\mathbf F}_{q^2}}

\newcommand\Ql{\bar{\mathbf Q}_l}
\newcommand\BQ{\mathbf Q}

\newcommand\BF{\mathbf F}

\newcommand\BZ{\mathbf Z}

\newcommand\BG{\mathbf G}

\newcommand\CE{\mathcal{E}}

\newcommand\CM{\mathcal{M}}

\newcommand\CO{\mathcal{O}}

\newcommand\CT{ \mathcal{T}}

\newcommand\Fu{\mathfrak u}
\newcommand\Fg{\mathfrak g}

\newcommand\iv{^{-1}}

\newcommand\wh{\widehat}
\newcommand\wt{\widetilde}
\newcommand\wg{^{\wedge}}

\newcommand\ol{\overline}
\newcommand\hra{\hookrightarrow}
\newcommand\lra{\leftrightarrow}

\newcommand\ssim{\! /\!\!\sim}
\newcommand\diag{\operatorname{diag}}
\newcommand\Id{\operatorname{Id}}

\newcommand\Ker{\operatorname{Ker}}
\newcommand\Hom{\operatorname{Hom}}
\newcommand\End{\operatorname{End}}

\newcommand\Ind{\operatorname{Ind}}

\newcommand\Res{\operatorname{Res}}
\newcommand\Irr{\operatorname{Irr}\,}

\newcommand\supp{\operatorname{supp}\,}
\newcommand\Lie{\operatorname{Lie}}
\newcommand\Tr{\operatorname{Tr}\,}

\newcommand\Ad{\operatorname{Ad}}
\newcommand\ad{\operatorname{ad}}

\newcommand\uni{_{\operatorname{uni}}}
\newcommand\nil{_{\operatorname{nil}}}

\newcommand\lp{\operatorname{\!\langle\!}}
\newcommand\rp{\operatorname{\!\rangle}}

\newcommand\nat{^{\natural}}

\newcommand\dc{\dot c}
\newcommand\ds{\dot s}
\newcommand\du{\dot u}
\newcommand\dx{\dot x}
\newcommand\dy{\dot y}
\newcommand\dz{\dot z}

\newcommand\hx{\hat x}
\newcommand{\isom}{\,\raise2pt\hbox{$\underrightarrow{\sim}$}\,}
\numberwithin{equation}{section}

\newtheorem{thm}{Theorem}[section]
\newtheorem{lem}[thm]{Lemma}
\newtheorem{cor}[thm]{Corollary}
\newtheorem{prop}[thm]{Proposition}

\def \para#1{\par\medskip\textbf{#1}
              \addtocounter{thm}{1}}

\def \remark#1{\par\medskip\noindent
                \textbf{Remark #1}
                \addtocounter{thm}{1}}

\begin{document}
\setlength{\baselineskip}{4.9mm}
\setlength{\abovedisplayskip}{4.5mm}
\setlength{\belowdisplayskip}{4.5mm}
\renewcommand{\theenumi}{\roman{enumi}}
\renewcommand{\labelenumi}{(\theenumi)}
\renewcommand{\thefootnote}{\fnsymbol{footnote}}

\parindent=20pt
\pagestyle{myheadings}
\markboth{SHOJI and SORLIN }{SYMMETRIC SPACE FOR $SL_n$}
\begin{center}
{\bf Subfield symmetric spaces for finite special linear groups} \\
\vspace{1cm}
Toshiaki Shoji and Karine Sorlin\footnote{
The second author would like to thank the JSPS for 
support which made this collaboration possible} \\ 
\vspace{0.5cm}
Graduate School of Mathematics \\
Nagoya University  \\
Chikusa-ku, Nagoya 464-8602,  Japan
\end{center}
\title{}
\maketitle
\begin{abstract}
Let $G$ be a connected algebraic group defined over a 
finite field $\Fq$.  For each irreducible character $\r$ of
$G(\BF_{q^r})$, we denote by $m_r(\r)$ the multiplicity of 
$1_{G(\Fq)}$ in the restriction of $\r$ to $G(\Fq)$. 
In the case where $G$ is reductive with connected center and 
is simple modulo center, Kawanaka determined $m_2(\r)$ for 
almost all cases, and then Lusztig gave a general 
formula for $m_2(\r)$.  
In the case where the center of $G$ is not connected, such a
result is not known.  In this paper we determine $m_2(\r)$, 
up to some minor ambiguity,  
in the case where $G$ is the special linear group.
\par
We also discuss, for any $r \ge 2$, the relationship
between $m_r(\r)$ with the theory of Shintani descent 
in the case where $G$ is a connected algebraic group. 
\end{abstract}
\par\medskip
\addtocounter{section}{-1}
\section{Introduction}
Let $G$ be a connected reductive group defined over a finite field 
$\Fq$ with Frobenius map $F$. We consider the finite group 
$G^{F^2}$ and 
its subgroup $G^F$.  The quotient space 
$G^{F^2}/G^F$ is regarded as an analogue of the symmetric space, and
is called the subfield symmetric space over a finite field.
The determination of spherical functions of $G^{F^2}/G^F$ is almost
equivalent to the determination of irreducible characters of the 
Hecke algebra $H(G^{F^2}, G^F)$.  For a class function  
$f$ on $G^{F^2}$, we denote by $m_2(f)$ the inner product of 
$\r$ with the induced character $\Ind_{G^F}^{G^{F^2}}1$.  The
classification of irreducible characters of $H(G^{F^2}, G^F)$ and 
the determination of their degrees are equivalent to the 
determination of $m_2(\r)$ for all irreducible characters $\r$ of 
$G^{F^2}$. 
\par
In [K2], Kawanaka computed $m_2(\r)$ in the case where $G$ 
is a classical group with connected center, or in the case where
$\r$ is unipotent and the characteristic is good.   
Extending Kawanaka's result, Lusztig gave in [L3] a closed formula
for $m_2(\r)$ valid for any $G$ which has the connected center and
is simple modulo its center. He expects that his formula is still 
valid for $G$ with disconnected center.  
In turn, Henderson studied in [H] the spherical functions of 
$G^{F^2}/G^F$ by making use of the theory of perverse sheaves, 
and described them in the case where $G = GL_n$, in which case 
$H(G^{F^2}, G^F)$ is abelian.    
\par
In this paper, we consider $G = SL_n$ with the standard
$\Fq$-structure, which is the first example of the 
disconnected center case.  Based on the parametrization of 
irreducible 
characters and the description of almost characters in [S3]
(which is valid under some restriction on $p$, for example, 
$p \ge n$), we determine
$m_2(\r)$ (Theorem 5.3) for any irreducible characters, 
up to some minor ambiguity.  Our result is consistent with 
Lusztig's conjectural formula modulo the ambiguity.  In particular, 
we have $m_2(\r) \in \{ 0,1, 2\}$.   
\par
Kawanaka's main idea for the computation of $m_2(\r)$, beside
the use of the results of Lusztig on $m_2(R_T(\th))$, is to connect 
it with the twisted Frobenius-Schur indicator through the twisting
operator.  In section 1, we generalize Kawanaka's result, and
discuss a connection of $m_2(\r)$ with Shintani descent.  This 
leads to a formula for $m_2(R_x)$ where $R_x$ is an almost 
character of $G^{F^2}$, which is regarded as a counter part
of Lusztig's formula for $m_2(\x_A)$ in [L3, 7], where $\x_A$ is the
characteristic function of character sheaves. 
In section 1, we also discuss a more general situation.
We define $m_r(\r)$ as the multiplicity of an irreducible character
$\r$ of $G^{F^r}$ with the induced character $\Ind_{G^F}^{G^{F^r}}1$
for any integer $ r \ge 2$.  We give some formula (Theorem 1.14) for 
$m_r(R_x)$ though it is not so effective as the $m_2$ case.
\par
The subsequent sections are devoted to the computation of 
$m_2(\r)$ for the case where $G = SL_n$.  We obtain the results
by applying the results in section 1, together with the computation of
$m_2(\wt\r|_{G^{F^2}})$ for irreducible characters $\wt\r$ of 
$GL_n(\Fqq)$.      
\par\bigskip\medskip
\begin{center}
{\sc Contents}
\end{center}
\par
1. $G(\Fq)$-invariants in $G(\Fqr)$-modules and Shintani descent.
\par
2. Parametrization of irreducible characters of $SL_n(\BF_{q^2})$.
\par
3. Almost characters of $SL_n(\BF_{q^2})$.
\par
4. Determination of $m_2(\r_{\ds, E}|_{G^{F^2}})$.
\par
5. Determination of $m_2(\r)$ for $\r \in \Irr SL_n(\BF_{q^2})$.

\par\medskip
\section{$G(\Fq)$-invariants in $G(\Fqr)$-modules and Shintani descent}
\para{1.1.}
For any finite group $\vG$ and an automorphism $F: \vG \to \vG$, we denote
by $\vG\ssim_F$ the set of $F$-twisted conjugacy classes in $\vG$, where
$x, y \in \vG$ are $F$-twisted conjugate if there exists $z \in \vG$ such
that $y = z\iv x F(z)$.
In the case where $F$ acts trivially on $\vG$, the set $\vG\ssim_F$ 
coincides with the set of conjugacy classes in $\vG$, which we denote by
$\vG\ssim$.
\par
For a connected algebraic group $X$ defined over $F_q$, and two
Frobenius maps $F_1, F_2$ on $X$ such that $F_1F_2 = F_2F_1$, 
we define a norm
map 
\begin{equation*}
N_{F_1/F_2} : X^{F_1}\ssim_{F_2} \to X^{F_2}\ssim_{F_1\iv}
\end{equation*}
as follows; 
for $x \in X^{F_1}$, we choose $\a \in X$ such that 
$x = \a\iv F_2(\a)$, and put $x' = F_1(\a)\a\iv$.  Then 
$x' \in X^{F_2}$ and the correspondence $x \to x'$ 
induces a bijective map 
$N_{F_1/F_2}$, which we call the norm map from
$X^{F_1}\ssim_{F_2}$ to $X^{F_2}\ssim_{F_1\iv}$.
\par
For a finite set $Y$, we denote by $C(Y)$ the $\Ql$-space
of all $\Ql$-valued functions on $Y$.
Then the norm map $N_{F_1/F_2}$ induces a linear isomorphism
\begin{equation*}
Sh_{F_1/F_2} = {N^{*-1}_{F_1/F_2}}: 
     C(X^{F_1}\ssim_{F_2}) \to C(X^{F_2}\ssim_{F_1\iv}),
\end{equation*}
which is called the Shintani descent from $X^{F_1}$ to $X^{F_2}$.
\para{1.2.}
Let $G$ be a connected algebraic group defined over a finite
field $\Fq$ with Frobenius map $F$.
We fix a positive integer $r$, and consider the group 
$H = G\times\cdots\times G$ ($r$-factors).
$H$ is endowed with the natural Frobenius map given by 
$(g_1, \dots, g_r) \mapsto (F(g_1), \dots, F(g_r))$, which we also
denote by $F$.  Let $F' = F\w : H \to H$ be a twisted Frobenius map 
on $H$, where 
$\w : H \to H, 
 (g_1, \dots, g_r) \mapsto (g_r, g_1, \dots, g_{r-1})$ 
is the cyclic permutation of factors.  
Since $\w^r = 1$ and $F\w = \w F$, we have $(F')^{rm} = F^{rm}$ 
for any $m \ge 1$.
\par
\begin{lem}
The map $G^{F^{rm}} \to H^{F^{rm}}, x \to (x, 1, \dots, 1)$
induces a bijection
\begin{equation*}
\tag{1.3.1}
f: G^{F^{rm}}\ssim_{F^r} \ \to H^{F^{rm}}\ssim_{F'}.
\end{equation*}
\end{lem}
\begin{proof}
Take $x = (x_1, \dots, x_r), y = (y_1, \dots, y_r) \in H^{F^{rm}}$.
If $x$ and $y$ are in the same class, there exists 
$z = (z_1, \dots, z_r)$ such
that $y_i = z_i\iv x_i F(z_{i-1})$ for $i \in \BZ/r\BZ$. 
Now assume that $x = (x_1, 1, \dots,1)$.  Then  
$z\iv xF'(z) = (y_1, 1, \dots,1)$ for $z \in G^{F^{rm}}$ 
if and only if $z = (z_1, F(z_1), \dots, F^{r-1}(z_1))$.
Moreover in this case, $y_1 = z_1\iv x_1F^r(z_1)$. 
This shows that the map $f$ is well-defined, and is injective. 
It is easy to see that each $F'$-conjugacy class in $H^{F^{rm}}$ 
contains a representative of the form $(x_1, 1, \dots, 1)$.
Hence $f$ is surjective.
\end{proof}
\para{1.4.}
For each $x \in G^{F^{rm}}$, $k \ge 1$, 
we put $N_k(x) = xF(x)\cdots F^{k-1}(x)$.
Then the map $G^{F^{rm}} \to G^{F^{rm}},  x \mapsto N_k(x)$ induces
a map $G^{F^{rm}}\ssim_F\ \to G^{F^{rm}}\ssim_{F^k}$, which we also 
denote by $N_k$. 
Let $\vD(H) \simeq G$ be the diagonal subgroup of $H$. The inclusion
$\vD(H)^{F^{rm}}\hra H^{F^{rm}}$ induces a map 
$d: \vD(H)^{F^{rm}}\ssim_F\ \to H^{F^{rm}}\ssim_{F'}$.
Then we have a commutative diagram
\begin{equation*}
\tag{1.4.1}
\begin{CD}
G^{F^{rm}}\ssim_{F^r} @>f>> H^{F^{rm}}\ssim_{F'} \\
@AN_rAA                                @AAdA  \\   
G^{F^{rm}}\ssim_F   @>f_0>>        \vD(H)^{F^{rm}}\ssim_F,
\end{CD}
\end{equation*}
where $f_0$ is the bijection induced from the isomorphism
$G \isom \vD(H)$.
This follows from the following relation for $x \in G^{F^{rm}}$, 
\begin{equation*}
(N_r(x), 1, \dots, 1) = y\iv(x, x, \dots, x)F'(y) 
\end{equation*}
with $y = (1, N_1(x), N_2(x), \dots, N_{r-1}(x))$.
\para{1.5.}
Concerning the norm maps, we have the following commutative diagram.
\begin{equation*}
\tag{1.5.1}
\begin{CD}
G^{F^{rm}}\ssim_{F^r} @> N_{F^{rm}/F^r}>> G^{F^r}\ssim \\
@AN_rAA                             @ AAjA  \\
G^{F^{rm}}\ssim_F @>N_{F^{rm}/F}>>        G^F\ssim, 
\end{CD}
\end{equation*}
where $j$ is the map induced from the inclusion $G^F \hra G^{F^r}$.
We show (1.5.1).  Let $\hx = G^{F^{rm}}$ and take $\a \in G$ such that
$\hx = \a\iv F(\a)$.  Then $N_{F^{rm}/F}(\hx)$ is represented by 
$x = F^{rm}(\a)\a\iv$. 
On the other hand, since $\hx' = N_r(\hx) = \a\iv F^r(\a)$,   
we see that $N_{F^{rm}/F^r}(\hx')$ is represented by 
$F^{rm}(\a)\a\iv$ which coincides with $j(x)$.  This shows 
the commutativity.
\para{1.6.}
Let  $\s' = F'|_{H^{F^{rm}}}$, and $\wt H^{F^{rm}}$ be 
the semidirect product of $H^{F^{rm}}$ with the cyclic group
$\lp\s'\rp$ of order $m$ generated by $\s'$. 
For a character $\x$ of $G^{F^{rm}}$, we define the character
$F(\x)$ by $F(\x)(F(g)) = \x(g)$, and similarly for $H$.
An irreducible character $\p$ of $H^{F^{rm}}$ is $F'$-stable
if and only if $\p$ is of the form that
\begin{equation*}
\tag{1.6.1}
\p = \x\otimes F(\x)\otimes\cdots\otimes F^{r-1}(\x)
\end{equation*}
for some
$F^r$-stable irreducible character $\x$ on $G^{F^{rm}}$. 
Let $V_i$ ($1 \le i \le r)$ be an irreducible $G^{F^{rm}}$-module 
for the irreducible character $F^{i-1}(\x)$.  Then there exists a 
linear isomorphism $T_i : V_i \to V_{i+1}$ such that 
$T_i\circ g = F(g)\circ T_i$ for any $g \in G^{F^{rm}}$ with 
$V_{r+1} = V_1$ and that $(T_rT_{r-1}\cdots T_1)^m = 1$.  
Let $\p$ be as in (1.6.1).  Then
$\p$ is afforded by
the $H^{F^{rm}}$-module $V_1\otimes V_2\otimes \cdots\otimes V_r$.
Let us define an action of $\s'$ on 
$V_1\otimes\cdots\otimes V_r$ by 
\begin{equation*}
\s' =  \w\circ (T_1\otimes T_2\otimes\cdots\otimes T_r),
\end{equation*}
where $\w$ is the cyclic permutation of factors given by 
\begin{equation*}
\w(x_1\otimes x_2\otimes\cdots\otimes x_r) = 
       x_r\otimes x_1\otimes\cdots\otimes x_{r-1}.
\end{equation*}
Then we have $\s'\circ h = F'(h)\circ \s'$ for $h \in H^{F^{rm}}$, 
and so 
$V_1\otimes\cdots\otimes V_r$ can be extended to an 
$\wt H^{F^{rm}}$-module.  We denote by $\wt\p$ the corresponding 
extension of $\p$ to $\wt H^{F^{rm}}$. 
\par
Let $\s = F|_{G^{F^{rm}}}$, and we consider $G^{F^{rm}}\lp\s\rp$ the 
semidirect product of $G^{F^{rm}}$ with the cyclic group $\lp\s\rp$ 
of order $rm$ generated by $\s$.
We define an action of $\s^r$ on $V_1$ by 
$\s^r = T_rT_{r-1}\cdots T_1$. 
Then $\s^r\circ g = F^r(g)\circ \s^r$ for any $g \in G^{F^{rm}}$, and 
the $G^{F^{rm}}$-module $V_1$ can be extended to a 
$G^{F^{rm}}\lp\s^r\rp$-module $\wt V_1$.
We denote by $\wt\x$ the corresponding extension 
of $\x$ to $G^{F^{rm}}\lp\s^r\rp$.  
We show the following lemma.
\begin{lem}  
Let $h = (g, 1, \dots, 1) \in H^{F^{rm}}$ with $g \in G^{F^{rm}}$.
Let $\x$ be an $F^r$-stable irreducible character of $G^{F^{rm}}$.
Then for 
$\p = \x\otimes F(\x)\otimes\cdots\otimes F^{r-1}(\x) 
      \in \Irr H^{F^{rm}}$, 
we have
\begin{equation*}
\wt\p(h\s') = \wt\x(g\s^r).
\end{equation*}  
\end{lem}
\begin{proof}
Let $v^{(1)}_1, \dots, v_n^{(1)}$ be a basis of $V_1$.  
We define a basis $v_1^{(i+1)}, \dots, v_n^{(i+1)}$ of $V_{i+1}$
inductively by $v_j^{(i+1)} = T_i(v_j^{(i)})$ for 
$i = 1,2, \dots r-1$.
Then we have
\begin{equation*} 
T_r(v_j^{(r)}) = T_r\cdots T_1(v_j^{(1)}) = \s^rv_j^{(1)}.
\end{equation*}
It follows that 
\begin{equation*}
h\s'\cdot v^{(1)}_{i_1}\otimes v^{(2)}_{i_2}\otimes\cdots\otimes 
          v^{(r)}_{i_r}
   = (g\s^rv^{(1)}_{i_r})\otimes v^{(2)}_{i_1}\otimes\cdots\otimes
           v^{(r)}_{i_{r-1}},                
\end{equation*}
and we have 
\begin{equation*}
\wt\p(h\s') = \Tr(h\s', V_1\otimes\cdots\otimes V_r) 
            = \Tr(g\s^r, V_1) = \wt\x(g\s^r). 
\end{equation*}
This proves the lemma.
\end{proof}
\para{1.8.}
Let $\x$ be an $F^r$-stable irreducible character of $G^{F^{rm}}$, 
and $\wt\x$ be its extension to 
$G^{F^{rm}}\lp\s^r\rp$ as in the 
previous lemma.  
Under the natural bijection 
$G^{F^{rm}}\ssim_{F^r} \simeq  G^{F^{rm}}\s^r\ssim$
via $x \lra x\s$, we have an isomorphism 
$C(G^{F^{rm}}\ssim_{F^r}) \simeq C(G^{F^{rm}}\s^r\ssim)$.
Thus $\wt\x|_{G^{F^{rm}}\s^r}$ defines an element in 
the space $C(G^{F^{rm}}\ssim_{F^r})$.  Put
\begin{equation*} 
R^{(m)}_{\wt\x} = Sh_{F^{rm}/F^r}(\wt\x|_{G^{F^{rm}}\s^r}).
\end{equation*}
Hence $R^{(m)}_{\wt\x}$ is a class function on $G^{F^r}$.
We have the following formula.
\begin{prop}
Under the notation as above, 
\begin{equation*}
\tag{1.9.1}
|G^{F^{rm}}|\iv\sum_{\hat g \in G^{F^{rm}}}\wt\x(N_r(\hat g)\s^r) = 
            |G^{F}|\iv\sum_{g\in G^{F}}R^{(m)}_{\wt\x}(g).
\end{equation*}
\end{prop}
\begin{proof}
Take $\hat g \in G^{F^{rm}}$.  Write $\hat g$ as 
$\hat g = \a\iv F(\a)$ and put $g = F^{rm}(\a)\a\iv$.  Then 
$g \in G^F$, and  we see that 
$\wt\x(N_r(\hat g)\s^r) = R_{\wt\x}(g)$ by (1.5.1).
Moreover, it is known that
\begin{equation*}
\sharp\{ x \in G^{F^{rm}} \mid x\iv \hat gF(x) = \hat g \}
   = \sharp\{ y \in G^F \mid y\iv gy = g\}.
\end{equation*}
The formula (1.9.1) is immediate from these two facts. 
\end{proof}
\para{1.10.}
Let $c_r^{(m)}(\wt\x)$ be the left hand side of (1.9.1), i.e., 
\begin{equation*}
\tag{1.10.1}
c_r^{(m)}(\wt\x) =  
 |G^{F^{rm}}|\iv\sum_{\hat g \in G^{F^{rm}}}\wt\x(N_r(\hat g)\s^r).
\end{equation*} 
Then $c_r^{(m)}(\wt\x)$ is a generalization of the twisted 
Frobenius-Schur indicator discussed in Kawanaka and Matsuyama [KM].  
In the case where $m = 1$, we simply write $c^{(1)}_r(\wt\x)$ as 
$c_r(\x)$.  Note that in this case, the extension does not enter the
formula, and we have 
\begin{equation*}
c_r(\x) = |G^{F^r}|\iv\sum_{g\in G^{F^r}}\x(N_r(g)).
\end{equation*}
If $r = 2$, $c_2(\x)$ coincides with the Frobenius-Schur indicator 
defined in [KM]. 
\par
Let us define, for a class function $f$ of $G^{F^r}$, 
\begin{equation*}
\tag{1.10.2}
m_r(f) = \lp f, \Ind_{G^F}^{G^{F^r}}\!1\rp = 
            |G^F|\iv\sum_{x \in G^F}f(x).
\end{equation*}
Then the identity (1.9.1) can be rewritten as 
\begin{equation*}
\tag{1.10.3}
c_r^{(m)}(\wt\x) = m_r(R_{\wt\x}^{(m)}).
\end{equation*}
We note that (1.10.3) is a generalization of the formula 
due to Kawanaka [K2, (1.1)].
 In fact, in the case where $m =1$, 
the Shintani descent $Sh_{F^r/F^r}$ coincides with the inverse 
of the twisting operator $t_1^*$ on $C(G^{F^r}\ssim)$ given 
in [K2], and so we have $R_{\x}^{(1)} = t_1^{*-1}\x$. 
Then (1.10.3) implies the following.
\begin{cor}  
Let the notations be as above.  
Then we have $ c_r(\x) = m_r(t_1^{*-1}\x)$.
\end{cor}
In the case where $r = 2$, this formula is 
nothing but the formula (1.1) in [K2]. 
\para{1.12.}
By Lemma 1.7, $\wt\x(N_r(\hat g)\s^r) = \wt\p(h\s')$ with
$h = (N_r(\hat g), 1, \dots,1) \in H^{F^{rm}}$.
As in 1.4, $h$ is $F'$-conjugate to 
$(\hat g, \dots, \hat g) \in \vD(H)^{F^{rm}}$, and so
\begin{equation*}
\wt\p(h\s') = \wt\p((\hat g, \dots, \hat g)\s').
\end{equation*}
On the other hand, under the isomorphism 
$\vD(H)^{F^{rm}} \simeq G^{F^{rm}}$, 
$V_1\otimes\cdots\otimes V_r$ is an $G^{F^{rm}}$-module, 
and its character $\x F(\x)\cdots F^{r-1}(\x)$
is $F$-stable. Moreover, we have $\s'\circ g = F(g)\circ \s'$
on $V_1\otimes\cdots\otimes V_r$ for any $g \in G^{F^{rm}}$.  
This implies that the action of 
$\s'$ defines a structure of  
$G^{F^{rm}}\lp\s\rp$-module  on $V_1\otimes\cdots\otimes V_r$, where
$\s$ acts by $\s'$ on it.  We denote the character of this module
by $\wt\p_0$, which is an extension of $\x F(\x)\cdots F^{r-1}(\x)$.  
Thus, we have
\begin{equation*}
\wt\p((\hat g, \dots, \hat g)\s') = \wt\p_0(\hat g\s).
\end{equation*} 
Now (1.9.1) can be rewritten as 
\begin{equation*}
\tag{1.12.1}
m_r(R_{\wt\x}^{(m)}) = 
|G^{F^{rm}}|\iv\sum_{\hat g \in G^{F^{rm}}}\wt\p_0(\hat g\s).
  \end{equation*}
\para{1.13.}
Let us define an inner product on $C(G^{F^{rm}}\ssim_F)$ by 
\begin{equation*}
\lp f, h \rp_{G^{F^{rm}}\s} = |G^{F^{rm}}|\iv\sum_{x \in G^{F^{rm}}}
                       f(x\s)\ol{h(x\s)}
\end{equation*}
for $f, h \in C(G^{F^{rm}}\ssim_F)$. 
Then the following orthogonality relations are known.
For any $F$-stable irreducible characters $\x, \x'$ of $G^{F^{rm}}$
and their extensions $\wt\x, \wt\x'$ to $G^{F^{rm}}\lp\s\rp$, 
\begin{equation*}
\tag{1.13.1}
\lp \wt\x, \wt\x'\rp_{G^{F^{rm}}\s} = \begin{cases}
           \th(\s) &\quad\text{ if } \wt\x' = \th\otimes \wt\x 
                     \text{ with } \th \in \Irr\lp \s\rp, \\
            0      &\quad\text{ if } \x' \ne \x.
                            \end{cases} 
\end{equation*}
Here in the left hand side, $\wt\x, \wt\x'$ are regarded as 
functions on $G^{F^{rm}}\s$ by restriction.
\par
For any $f \in C(G^{F^{k}}\ssim_F)$,  we put
\begin{equation*}
\wt M_{k}(f) = \lp f, \wt 1\rp\,_{G^{F^{k}}\s} 
              = |G^{F^{k}}|\iv\sum_{x \in G^{F^{k}}}f(x\s),
\end{equation*}
where $\wt 1$ means the restriction of 
the unit character of $G^{F^{k}}\lp\s\rp$ to $G^{F^{k}}\s$.
We also put, for a class function $h$ of $G^{F^k}$, 
\begin{equation*}
M_{k}(h) = \lp h, 1\rp\,_{G^{F^{k}}} 
              = |G^{F^{k}}|\iv\sum_{x \in G^{F^{k}}}h(x).
\end{equation*}
The following statement is immediate from (1.13.1).
\par\medskip\noindent
(1.13.2) \ Let $\r$ be an $F$-stable character of $G^{F^{k}}$, 
and $\wt\r$ its extension to $\wt G^{F^{k}}$.  Then 
we have $|\wt M_{k}(\wt\r)| \le M_{k}(\r)$.  
Moreover, if $M_{k}(\r) = 1$, then $\wt M_{k}(\wt\r)$ is a
$k$-th root of unity.
\par\medskip
We have the following theorem.
\begin{thm} 
Let $\x$ be an $F^r$-stable irreducible character of $G^{F^{rm}}$, and 
$\wt\x$ an extension of $\x$ to $G^{F^{rm}}\lp\s^r\rp$.
Let $\wt\p_0$ be the extension of $\x F(\x)\cdots F^{r-1}(\x)$ to
$G^{F^{rm}}\lp\s\rp$ as in 1.12.  
Put $Sh_{F^{rm}/F^r}(\wt\x|_{G^{F^{rm}}\s^r}) = R^{(m)}_{\wt\x}$. 
\begin{enumerate}
\item
We have 
$c_r^{(m)}(\wt\x) = m_r(R^{(m)}_{\wt\x}) = \wt M_{rm}(\wt\p_1)$.
In particular, 
\begin{equation*}
|m_r(R^{(m)}_{\wt\x})| \le M_{rm}(\x F(\x)\cdots F^{r-1}(\x)).
\end{equation*}
Furthermore, if $M_{rm}(\x F(\x)\cdots F^{r-1}(\x)) = 1$, we have
$|m_r(R^{(m)}_{\wt\x})| = 1$.
 
\item
Assume that $r = 2$.  Then there exists a $2m$-th root of unity $\z$
such that 
\begin{equation*}
m_2(R^{(m)}_{\wt\x}) = \begin{cases}
                            \z &\quad\text{ if } \ol\x = F(\x), \\
                             0 &\quad\text{ otherwise}, 
                       \end{cases}
\end{equation*}
where $\ol\x$ is the complex conjugate of the character $\x$.
\end{enumerate}
\end{thm}
\begin{proof}
The equality $m_r(R^{(m)}_{\wt\x}) = \wt M_{rm}(\wt\p_0)$ in (i) 
follows from (1.12.1).  The inequality in (i) follows from (1.13.2).
Assume that $r = 2$.  Then we have
\begin{equation*}
M_{2m}(\x F(\x)) = \lp \x F(\x), 1\rp_{G^{F^{2m}}} = 
                   \lp F(\x), \ol\x\rp_{G^{F^{2m}}}
                 = \begin{cases}
                     1 &\quad\text{ if } F(\x) = \ol\x, \\
                     0 & \quad\text{ otherwise}.
                   \end{cases}. 
\end{equation*}
So the assertion (ii) follows from (1.13.2).  This proves the
theorem.
\end{proof}
\para{1.15.}
In the case where $r = 2$, we determine the quantity 
$\z = c_2^{(m)}(\wt\x)$ more explicitly. 
Let $\x$ be an $F^2$-stable irreducible character of $G^{F^{2m}}$ 
and $\wt\x$ its extension to $G^{F^{2m}}\lp\s^2\rp$ 
as in the theorem.  Let us assume that $F(\x) = \ol\x$.
We follow the setting in 1.6.  In particular $V_1$ (resp. $V_2$)
is a $G^{F^{2m}}$-module affording $\x$ (resp. $F(\x)$).
Since $F(\x) = \ol\x$, the subspace $W = (V_1\otimes V_2)^{G^{F^{2m}}}$  
of $G^{F^{2m}}$-invariant vectors in $V_1\otimes V_2$ is of dimension 1.
The map $\s': V_1\otimes V_2 \to V_1\otimes V_2, 
v_1\otimes v_2 \mapsto T_2(v_2)\otimes T_1(v_1)$ preserves the space 
$W$, and the eigenvalue of $\s'$ on $W$ coincides with 
$\z = c_2^{(m)}(\wt\x)$.  
The map $\s^2 = T_2T_1 : V_1 \to V_1$ extends the $G^{F^{2m}}$-module
$V_1$ to the $G^{F^{2m}}\lp\s^2\rp$-module $\wt V_1$ 
affording the character $\wt\x$.
\par
The $G^{F^{2m}}$-module $V_2$ can be identified with 
$V_1$ by replacing the action of $g \in G^{F^{2m}}$ by $F(g)$.
Under this identification, we may take $T_1 = \Id_{V_1}$ 
and $T_2 = \s^2$ on $V_1$. Hence we have
$\s'(v_1\otimes v_2) = \s^2(v_2)\otimes v_1$. 
Now the averaging operator 
$V_1\otimes V_2 \to W, v \mapsto 
  |G^{F^{2m}}|\iv\sum_{g \in G^{F^{2m}}}g\cdot v$
determines a bilinear form $B: V_1 \times V_1 \to \Ql$ (up to scalar)
having the following properties.
\begin{align*}
\tag{1.15.1}
B(g\cdot v_1, F(g)\cdot v_2) &= B(v_1, v_2) \quad\text{ for } 
                 g \in G^{F^{2m}}, v_1, v_2 \in V_1 \\ 
B(\s^2(v_2), v_1) &= \z B(v_1, v_2) \quad\text{ for } v_1, v_2 \in V_1.
\end{align*}
Conversely, if there exists such a bilinear from on $V_1$, 
this form coincides with $B$ up to scalar.  Hence $\z$ determines the
value $c_2^{(m)}(\wt\x)$.
\par
The extension $\wt\x$ of $\x$ is determined by the choice of $T_1, T_2$
such that $(T_2T_1)^m = \Id_{V_1}$.
If we replace $T_1$ by a scalar multiple $\xi T_1$
for an $m$-th root of unity $\xi$, it gives a different extension of 
$\wt\x'$ of $\x$. By changing $\wt\x$ by $\wt\x'$, the eigenvalue $\z$
of $\s'$ on $W$ is replaced by $\xi\z$.
Summing up the above arguments, we have the following refinement of 
Theorem 1.14, which is a generalization of Theorem 2.1.3 in [K2].
\begin{cor} 
Let $\x$ be an $F^2$-stable irreducible character of $G^{F^{2m}}$
and $\wt\x$ an extension of $\x$ to $G^{F^{2m}}\lp\s^2\rp$.
\begin{enumerate}
\item We have
\begin{equation*}
c_2^{(m)}(\wt\x) = \begin{cases}
                       \z &\quad\text{ if } F(\x) = \ol\x, \\
                       0  &\quad\text{ otherwise}, 
                   \end{cases}
\end{equation*}      
where $\z$ is an $2m$-th root of unity.
\item
Assume that $F(\x) = \ol\x$.   
Let $\z_0$ be a primitive $2m$-th root of unity in $\Ql$.
Then there exists a unique extension $\wt\x$ of $\x$ such that
$c_2^{(m)}(\wt\x) = 1$ or $\z_0$.  Let $V_1$ be the 
$G^{F^{2m}}\lp\s^2\rp$-module affording $\wt\x$.
Then $c_2^{(m)}(\wt\x) = 1$ (resp $\z_0)$ if and only if there exists 
a non-zero bilinear form $B(\cdot, \cdot)$ on $V_1$ satisfying
(1.15.1) with $\z = 1$ (resp. $\z = \z_0$).
\end{enumerate}
\end{cor} 
\para{1.17.}
In the case where $G$ is a connected reductive group with connected 
center, Lusztig defined in [L1] almost characters of $G^F$.  
In the case where $G$ is a special linear group $SL_n$ with $F$ of 
split type,  almost characters are also formulated in [S3].  
In either case, the set of almost characters coincides with 
the set of $Sh_{F^m/F}(\wt\x|_{G^{F^m}\s})$, up to an $m$-th root 
of unity multiple, 
for sufficiently divisible $m$, where $\x$ runs over all the 
$F$-stable irreducible characters of $G^{F^m}$.  
We denote by $R_{\x}$ the almost character of $G^F$ corresponding to 
$\x$.  As a corollary to Theorem 1.14, we have the following result.
\begin{cor} 
Assume that $G$ is either a connected reductive group with connected center, 
or $SL_n$ with $F$ of split type.  Let $R_{\x}$ be the almost character
of $G^{F^2}$ associated to an $F^2$-stable irreducible character 
$\x$ of $G^{F^{2m}}$.  Then we have
\begin{equation*}
\tag{1.18.1}
m_2(R_{\x}) = \begin{cases}
                \z &\quad\text{ if } F(\x) = \ol\x, \\
                 0 &\quad\text{ otherwise},
               \end{cases}
\end{equation*} 
where $\z$ is a certain $2m$-th root of unity.
\end{cor} 
\remark{1.19.}
In [L3, Prop. 7.2], Lusztig proved a formula concerning the characteristic 
functions 
of character sheaves as follows.  Let $A$ be an $F^2$-stable character
sheaf of a connected reductive group $G$.  We denote by 
$\x_{A, \f_A} \in C(G^{F^2}\ssim)$ the characteristic 
function of $A$ with respect to an isomorphism 
$\f_A: (F^2)^*A \isom A$.  Then under the assumption that $q$ is
sufficiently large (and that $\x_{A, \f_A}$ can be written 
as a linear combination of irreducible characters with 
cyclotomic integers coefficients), there exists a choice of $\f_A$
such that
\begin{equation*}
\tag{1.19.1}
m_2(\x_{A, \f_A}) = \begin{cases}
                        (-1)^{\dim\supp A} &\quad\text{ if } 
                                  F^*(A) \simeq DA, \\
                         0                 &\quad\text{ otherwise,}
                     \end{cases} 
\end{equation*} 
where $DA$ is the Verdier dual of $A$.
Since the proof depends on the asymptotic behavior of $q \to \infty$, 
the condition on $q$ is considerably large. In the case where $G$ has
a connected center, using the description of $m_2(\x)$ for any 
irreducible character $\x$ of $G^{F^2}$ in [L3], (1.18.1) can be
verified directly.  In [S2], it was shown that almost characters 
coincide with the characteristic functions of character
sheaves whenever $G$ has a connected center.  A similar result was
also shown in [S4] for $SL_n$ with $F$ of split type.  Hence 
the formula (1.18.1) is a counter part of (1.19.1) to 
almost characters, which works without any assumption on $q$.
Also, Theorem 1.14 (ii) is regarded as an extension of 
(1.19.1) to arbitrary connected algebraic groups.
\para{1.20.}
As a special case of the situation discussed in Theorem 1.14 (i), 
we consider the case where $G = GL_n$ with the standard or
non-standard Frobenius map $F$ over $\Fq$.  
Irreducible characters of $G^{F^r}$ is described 
as follows.  Let $G^* \simeq GL_n$ be the dual group of $G$.
For each $F^r$-stable semisimple class $\{ s\}$, choose 
a representative $s \in G^{*F^r}$.  Let $T^*$ be a maximally
split maximal torus in $Z_{G^*}(s)$.
Let $W = N_{G^*}(T^*)/T^*$ be the Weyl group of $G^*$, and put
$W_s = \{ w \in W \mid w(s) = s \}$.  Then $W_s$ is the Weyl group
of $Z_{G^*}(s)$, and $F^r$ acts naturally on $W_s$, which we denote 
by $\d$.
Let $(\Irr W_s)^{\d}$ be the set of $F^r$-stable irreducible
representations of $W_s$.  For each $E \in (\Irr W_s)^{\d}$, 
we fix an extension $\wt E$ of $E$ to the semidirect group 
$W_s\lp\d\rp$, where $\lp\d\rp$ is the infinite cyclic 
group with generator $\d$.   Put 
\begin{equation*}
R_{s, \wt E} = 
|W_s|\iv\sum_{w \in W_s}\Tr(w\d, \wt E)R_{T^*_w}(s),
\end{equation*} 
where $R_{T^*_w}(s)$ denotes  the Deligne-Lusztig character
$R_{T_w}(\th)$ under the natural correspondence 
$(s, T_w^*) \lra (\th, T_w)$.
\par
It is known, under a suitable choice of the extension, 
$\pm R_{s, \wt E}$ gives rise to an irreducible character of 
$G^{F^r}$, which we denote by $\r_{s,E}$. 
Then the set $\Irr G^{F^r}$ of irreducible characters of $G^{F^r}$ is 
given as 
\begin{equation*}
\Irr G^{F^r} = 
  \coprod_{\{s\}} \{ \r_{s, E} \mid E \in (\Irr W_s)^{\d}\},
\end{equation*}
where $\{s\}$ runs over $F^r$-stable semisimple conjugacy classes 
in $G^*$.  
\par
Let $(s, T^*)$ be as above.  We choose an $F^r$-stable maximal torus
of $G^{F^r}$ which is dual to $T^*$, and let $B$ be a Borel subgroup 
of $G$ containing $T$.  We choose an integer $m > 0$ such that 
$F^{mr}$ leaves $B$ invariant. 
One can find a linear character $\th$ of $T^{F^{rm}}$ corresponding 
to $s \in T^{*F^{rm}}$.
Then we have 
\begin{equation*}
\End_{G^{F^{rm}}}\bigl(\Ind_{B^{F^{rm}}}^{G^{F^{rm}}}\wt\th\bigr) 
                      \simeq \Ql[W_s], 
\end{equation*} 
where $\wt\th$ is the lift of $\th$ to the linear character of 
$B^{F^{rm}}$. 
Let us denote by $\x_{\th, E}$ the irreducible constituent of 
$\Ind_{B^{F^{rm}}}^{G^{F^{rm}}}\wt\th$ corresponding to 
$E \in \Irr W_s$.
Then $\x_{\th, E}$ is $F^r$-stable if and only if 
$E \in (\Irr W_s)^{\d}$,
and in which case, 
$Sh_{F^{rm}/F^r}(\wt\x_{\th, E}|_{G^{F^{rm}}\s^r})$ 
coincides with 
$\r_{s, E}$ up to a scalar multiple. 
Thus under this setting, Theorem 1.14 (i) can be rewritten as
follows.
\begin{cor} 
Let $G = GL_n$ with the standard or non-standard Frobenius map $F$.  
Then for each $\r_{s, E} \in \Irr G^{F^r}$, we have
\begin{equation*}
m_r(\r_{s, E}) \le 
      M_{rm}(\x_{\th, E}F(\x_{\th,E})\cdots F^{r-1}(\x_{\th,E})).
\end{equation*} 
Moreover, if 
$M_{rm}(\x_{\th,E}F(\x_{\th,E})\cdots F^{r-1}(\x_{\th,E})) = 1$, 
we have $m_r(\r_{s,E}) = 1$.
\end{cor}
\par\medskip
\section{Parametrization of irreducible characters of $SL_n(\Fqq)$}
\para{2.1.}
In the remainder of this paper, we assume that $\wt G = GL_n$ 
and $G = SL_n$ with
Frobenius maps  $F$ with respect to the standard $\Fq$-structures.
We assume that $p$ is large enough so that the results in [S3] can 
be applied.  For example $p \ge n$ is enough in our case.
Let $\wt G^*$ (resp. $G^*$) be the dual group of $\wt G$ (resp. $G$).
Then $\wt G^* \simeq GL_n$, and $G^* \simeq \wt G^*/\wt Z^*$, where 
$\wt Z^*$ is the center of $\wt G^*$.  The inclusion map 
$G \hra \wt G$ induces a natural surjection 
$\pi : \wt G^* \to G^* $.
  As in the case of $\wt G$, the set $\Irr G^{F^2}$
is partitioned as
\begin{equation*}
\Irr G^{F^2} = \coprod_{\{ s\}}\CE(G^{F^2}, \{  s\}),
\end{equation*}
where $\{ s\}$ runs over $F^2$-stable semisimple classes in $G^*$.
Take $s$ such that $F^2(s) = s$. 
Let $T^*$ be an $F^2$-stable maximal torus of $Z_{G^*}(s)$ such that
$T^*$ is contained in an $F^2$-stable Borel subgroup of $Z_{G^*}(s)$.
Let $\wt T^*$ be an $F^2$-stable maximal torus of $\wt G^*$ such that 
$\pi(\wt T^*) = T^*$. 
Then $W = N_{\wt G^*}(\wt T^*)/\wt T^*$ is naturally identified with
$N_{G^*}(T^*)/T^*$.
Put 
\begin{equation*}
W_s = N_{Z_{G^*}(s)}(T^*)/T^*, \quad 
W_s^0 = N_{Z^0_{G^*}(s)}(T^*)/T^*. 
\end{equation*}
 Then $W_s^0$ is the Weyl group
of $Z^0_{G^*}(s)$.  
Now $W_s$ can be decomposed as
$W_s \simeq W^0_s \rtimes \Om_s$, where 
$\Om_s = Z_{G^*}(s)/Z^0_{G^*}(s)$ is a cyclic group.  
If we choose $\ds \in \wt T^*$ such that 
$\pi(\ds) = s$, then
$W_s^0$ is naturally identified with 
$W_{\ds} = \{ w \in W \mid w(\ds) = \ds \}$.
\par
$F^2$ acts naturally on $W_s$.  We denote by $\d$ this action and 
consider the semidirect product $W_s\lp\d\rp$, where 
$\d w\d\iv = F^2(w)$.  
$\d$ stabilizes $W_s^0$ and $\Om_s$.
\para{2.2.}
For each $E \in \Irr W_s^0$, let $\Om_s(E)$ be the stabilizer of 
$E$ in $\Om_s$. 
Assume that the $\Om_s$-orbit of $E$ is $\d$-stable.  Put
\begin{equation*}
\wt\Om_s(E) = \{ u \in \Om_s \mid {}^{u\d}E = E \}.
\end{equation*}
Then one can write $\wt\Om_s(E) = \Om_s(E)a$ for some $a \in \Om_s$.
Since $\Om_s(E)$ is abelian, $\Om_s(E)$ is $\d$-stable, and 
$\Om_s(E)$ acts on $\wt\Om_s(E)$ by $(z, u) \mapsto z\iv u\d(z)$ 
for $z \in \Om_s(E)$ and $u \in \wt\Om_s(E)$. 
We denote by $\wt\Om_s(E)_{\d}$ the set of equivalent classes under 
this action.
It is easy to see that $\wt\Om_s(E)_{\d}$ can be identified with
the set $\Om_s(E)_{\d}a$, where $\Om_s(E)_{\d}$ is the 
largest quotient of $\Om_s(E)$ on which $\d$ acts trivially. 
Let $\ol{\Irr} W^0_s$ be the set of $\Om_s$-orbits in the set 
$\Irr W^0_s$.  We denote by 
$(\ol{\Irr} W^0_s)^{\d}$ the set of $\d$-stable orbits in $\Irr W^0_s$.     
\par
For each pair $(s, E)$ with $E \in (\ol{\Irr}W^0_s)^{\d}$, put
\begin{equation*}
\ol\CM_{s, E} = \Om_s^{\d}(E)\wg \times \wt\Om_s(E)_{\d}, 
\end{equation*}
where $\Om_s^{\d} = \{ u \in \Om_s, \d(u) = u \}$ and  
$\Om_s^{\d}(E)$ is the stabilizer of $E$ in $\Om_s^{\d}$, 
and $\Om_s^{\d}(E)\wg$ is the set of irreducible characters of 
$\Om_s^{\d}(E)$.
It is known by [S3] that there exists a natural bijection 
\begin{equation*}
\tag{2.2.1}
\CE(G^{F^2}, \{ s\}) = 
     \coprod_{E \in (\ol{\Irr} W^0_s)^{\d}} \ol\CM_{s, E}
\end{equation*}
We denote by $\r_{\e, z}$ the irreducible character of $G^{F^2}$
corresponding to $(\e, z) \in \ol\CM_{s, E}$.
\par
The above parametrization satisfies the following properties;
The set of $G^{*F^2}$-conjugacy classes in  the set $\{ s\}^{F^2}$
is in bijection with $(\Om_s)_{\d}$.  
For each $x \in (\Om_s)_{\d}$, take a representative 
$\dx \in \Om_s$, choose $g_x \in G^*$ such that 
$g_x\iv F^2(g_x) = \dx$,  and
 put $s_x = g_xsg_x\iv$.  Then $F^2(s_x) = s_x$, and 
$g_xT^*g_x\iv = T^*_x$ is a maximally split torus in $Z_{G^*}(s_x)$.
We define $W_{s_x}^0$ in a similar way as $W_s^0$.
Under the isomorphism $W^0_s \to W^0_{s_x}$ 
induced by $\ad g_x$, the action of $x\d$ on $W_s^0$ is 
transferred to the action of $F^2$ on $W^0_{s_x}$. 
Hence each $x\d$-stable irreducible character $E'$ of $W^0_s$ 
determines the $F^2$-stable irreducible character $E''$ of
$W^0_{s_x}$.
Take an $F^2$-stable element $\ds_x$ such that 
$\pi(\ds_x) = s_x$.  We consider the irreducible 
character $\r_{\ds_x, E''}$ of $G^F{^2}$ as in 1.19, 
which we denote by $\r_{\ds_x, E'}$, by abuse of the 
notation. 
\par
It is known by [S3, (4.4.2)] that there exists a natural 
bijection 
\begin{equation*}
\tag{2.2.2}
f: \coprod_{E \in (\ol{\Irr} W^0_s)^{\d}} \wt\Om_s(E)_{\d}
   \simeq \coprod_{ x\in (\Om_s)_{\d}}(\Irr W^0_s)^{\dx\d}/\Om_s^{\d},
\end{equation*} 
where in the right hand side, $(\Irr W^0_s)^{\dx\d}/\Om_s^{\d}$
means the set of $\Om_s^{\d}$-orbit of $\dx\d$-stable irreducible
characters of $W_s^0$.
The bijection is described as follows.  Take $E$ in a $\d$-stable
$\Om_s$-orbit in $\Irr W^0_s$.  For each $\dy \in \wt\Om_s(E)$, there 
exists $\dx \in \Om_s$ and $z \in \Om_s$ such that
$\dy = z\iv \dx\d(z)$. 
Then $E_x = {}^zE \in (\Irr W^0_s)^{\dx\d}$.  The correspondence
$(E, y) \mapsto (x, E_x)$ gives rise to the required 
bijection $f$.
\par
Under the above setting, we have
\begin{equation*}
\tag{2.2.3}
\r_{\ds_x, E_x}|_{G^{F^2}} = \sum_{\e \in \Om_s^{\d}(E)\wg}\r_{\e,y}. 
\end{equation*}
\par
Let $\CT_{s_x, E_x}$ be the set of irreducible characters occurring 
in the restriction of $\r_{\ds_x, E_x}$ to $G^{F^2}$.  We also denote
by $\ol\CT_{s,E}$ the set of $\r_{\e,y}$ for 
$(\e, y) \in \ol\CM_{s, E}$.  Then (2.2.2) implies that 
\begin{equation*}
\ol\CT_{s, E} = \coprod_{(x, E_x)} \CT_{s_x, E_x},
\end{equation*}
where $(x, E_x)$ runs over all the pairs corresponding to 
$(y, E)$ with $y \in \wt\Om_s(E)_{\d}$ under the map $f$.
\remark{2.3.} \
In [S3, 4.5], the parameter set $\ol\CM_{s, E}$ is defined as 
$\Om_s^{\d}(E)\wg \times \Om_s(E)_{\d}$.  Since 
$\wt\Om_s(E)_{\d} = \Om_s(E)_{\d}a$,  this set is in bijection with
$\ol\CM_{s, E}$ in this paper.  However, the bijection 
depends on the choice of $a \in \Om_s$, and the definition of 
$\ol\CM_{s, E}$ in this paper is more convenient for later 
applications.  
\para{2.4.}
We describe the decomposition of $\r_{\ds_x, E_x}|_{G^{F^2}}$ 
in (2.2.3) more precisely.
It is known by [L2] that $\CT_{s_x,E_x}$ is in bijective 
correspondence with  $\Om_s^{\d}(E)\wg$. 
This bijection is given as follows.    
The abelian group  $\wt G^{F^2}/G^{F^2}$ acts transitively on
$\CT_{s_x, E_x}$ by the conjugation action.  
Also its dual group $(\wt G^{F^2}/G^{F^2})\wg$ acts on 
$\Irr \wt G^{F^2}$ by $(\th, \wt\r) \mapsto \th\otimes\wt\r$ for a linear
character $\th \in (\wt G^{F^2}/G^{F^2})\wg$ and 
$\wt\r \in \Irr \wt G^{F^2}$.
Then for $\r_0 \in \CT_{s_x, E_x}$, the stabilizer of $\r_0$ in 
$\wt G^{F^2}/G^{F^2}$ and the stabilizer of $\r_{\ds_x, E_x}$ in 
$(\wt G^{F^2}/G^{F^2})\wg$ are orthogonal to each other under the
natural duality pairing 
$\wt G^{F^2}/G^{F^2} \times (\wt G^{F^2}/G^{F^2})\wg \to \Ql$ 
(cf. [L2, 9]).
Let $I(\r_{\ds_x, E_x})$ be the stabilizer of $\r_{\ds_x, E_x}$ in 
$(\wt G^{F^2}/G^{F^2})\wg$.  Then, under the choice of $\r_0$, 
the set $\CT_{s_x,E_x}$ is in natural bijection with 
$I(\r_{\ds_x, E_x})\wg$. 
\par
We show that $I(\r_{\ds_x,E_x})$ is isomorphic to $\Om_s^{\d}(E)$.
First note that there exists a natural isomorphism 
\begin{equation*}
\tag{2.4.1}
\wt Z^{*F^2} \simeq \Hom (\wt G^{F^2}/G^{F^2},\Ql^*) = 
         (\wt G^{F^2}/G^{F^2})\wg.
\end{equation*}
If $z$ is an element in $\wt Z^{*F^2}$ corresponding to 
$\th \in (\wt G^{F^2}/G^{F^2})\wg$ under the above isomorphism,
then $\th$ maps $\CE(\wt G^{F^2}, \{\ds_x\})$ onto 
$\CE(\wt G^{F^2}, \{ z\ds_x\})$.
Put
\begin{equation*}
\wt Z^{*F^2}_{s_x} = \{ z \in \wt Z^{*F^2} \mid 
    \text{ $z\ds_x$ is conjugate to $\ds_x$ under $\wt G^*$} \}, 
\end{equation*}
which does not depend on the choice of $\ds_x$ for $s_x$.
Then, under the identification in (2.4.1),  
$I(\r_{\ds_x,E_x})$ is regarded as a subgroup of $\wt Z^{*F^2}_{s_x}$.
Here we have a natural isomorphism
\begin{equation*}
\tag{2.4.2}
\begin{CD}
\w_{s_x}: \Om^{\d}_s = \Om^{x\d}_s  @>\ad g_x>> 
       Z_{G^*}(s_x)^{F^2}/Z^0_{G^*}(s_x)^{F^2} @>>> 
        \wt Z_{s_x}^{*F^2} 
\end{CD}
\end{equation*}
defined as follows. For each $z \in Z_{G^*}(s)^{xF^2}$, choose
$\dz \in \wt G^{*xF^2}$ such that $\pi(\dz) = z$. 
Then ${}^{g_x}(\ds\iv\dz\ds\dz\iv) \in \wt Z_{s_x}^{*F^2}$, 
and this induces the required isomorphism since
$\pi(Z_{\wt G^*}(\ds_x)) = Z^0_{G^*}(s_x)$.
Now under the identification in  (2.4.1), (2.4.2), we may see that
$I(\r_{\ds_x,E_x})$ is a subgroup of $\Om_s^{\d}$, and in fact, 
$I(\r_{\ds_x, E_x})$ coincides with the stabilizer of $E_x$ in
$\Om_s^{\d}$.  Thus we have 
$I(\r_{\ds_x, E_x}) = \Om_s^{\d}(E_x) = \Om_s^{\d}(E)$.
\para{2.5.}
The bijection between $\CT_{s_x, E_x}$ and $\Om_s^{\d}(E)\wg$ given in 
2.4 depends on the choice of $\r_0 \in \CT_{s_x, E_x}$.
We have to choose a specific $\r_0$ for each $\CT_{s_x, E_x}$.
This problem is reduced to a certain special case, and is solved
by the aide of generalized Gelfand-Graev characters.
\par
Let $\Fg$ be the Lie algebra of $G$ with Frobenius map $F$.
We have a bijection $\log: G\uni \to \Fg\nil$ by $v \mapsto v-1$, 
where $G\uni$ (resp. $\Fg\nil$) is the unipotent variety of $G$
(resp. nilpotent variety of $\Fg)$.
Let $N$ be a nilpotent element in $\Fg^F$.
By Dynkin-Kostant theory, there exists a natural grading 
$\Fg = \bigoplus_{i \in \BZ}\Fg_i$ associated to $N$. 
Let $\Fu_i = \bigoplus_{j \ge i}\Fg_j$.
Then one can find an $F$-stable parabolic subgroup $P = LU_1$ 
associated to $N$, where $L$ is an $F$-stable Levi subgroup of $P$ 
with $\Lie L = \Fg_0$,  and 
$U_1$ is the unipotent radical of $P$ with $\Lie U_1 = \Fu_1$.  
Moreover we have $N \in \Fg_2$. 
Let $k$ be an algebraic closure of $\Fq$.
According to Kawanaka [K1], 
there exists an $F$-stable subspace $\Fu$ ($\Fu_{1.5}$ in the notation 
of [S3]) 
of $\Fu_1$ containing $\Fu_2$ and an $F$-equivariant linear map 
$\la: \Fu \to k$ satisfying the following. 
There exists an $F$-stable connected unipotent subgroup $U$ of 
$U_1$ such that $\log(U) = \Fu$ and that 
the map $\la\circ\log : U \to k$ turns out to be an 
$F$-stable homomorphism of $U$.  
We define a linear character $\vL_N$ of $U^{F^2}$ by 
$\vL_N = \p_2\circ\la\circ\log$, where $\p_2: \Fqq \to \Ql^*$ is the 
additive character defined by $\p_2 = \p\circ\Tr_{\Fqq/\Fq}$ for a 
non-trivial additive character $\p: \Fq \to \Ql^*$.
The generalized Gelfand-Graev character $\vG_N$ of $G^{F^2}$ associated
to $N$ is defined as 
\begin{equation*}
\vG_N = \Ind_{U^{F^2}}^{G^{F^2}}\vL_N.
\end{equation*}
The character $\vG_N$ depends only
on the $G^{F^2}$-conjugacy class of $N$.
\par
We now consider the following special setting for the set 
$\ol\CM_{s, E}$ determined by the pair $(s, E)$. 
\par\medskip\noindent
(2.5.1) \ $W^0_s$ is isomorphic to $S_b \times \cdots\times  S_b$ ($t$-times)
with $b = n/t$, and $\Om_s \simeq \lp w_0\rp$, where $w_0$ is an
element of order $t$ in $W_{s}$ permuting the factors of 
$W^0_s$ transitively. 
Moreover, $E \in \Irr W^0_s$ is of the form 
\begin{equation*}
E = E_1\boxtimes E_1\boxtimes\cdots\boxtimes E_1 \quad\text{with}
      \quad  E_1 \in \Irr S_b. 
\end{equation*}
Then $E$ is $\Om_s$-stable, and we have $\Om_s = \Om_s(E)$.
Now it is known that there exists a unique irreducible character
$\r_0$ such that $\r_0$ occurs both in $\vG_N$ and in 
$\r_{\ds_x, E}|_{G^{F^2}}$.  By using this $\r_0$, one obtains a
bijection $\CT_{s_x, E_x} \lra  \Om_s^{\d}(E)\wg$ as in 2.4.  This is
the parametrization given in (2.2.3), where if $(x, E_x)$ corresponds
to $(E,y)$ by (2.2.2), then $\r_{\e, y}$ corresponds to 
$\e \in \Om_s^{\d}(E)\wg$.  
\par
By the arguments in [S3, 4.5], the parametrization of $\CT_{s_x, E_x}$
in the general case is reduced to the case given in (2.5.1).
Accordingly, $\r_0$ is determined for each $\CT_{s_x, E_x}$.
However, note that this parametrization still depends on the choice of 
a nilpotent element $N$ in $\Fg$.  In what follows, we assume that
\par\medskip\noindent
(2.5.2) \ Each nilpotent element $N \in \Fg^{F}$ 
is taken to be a Jordan normal form.
\para{2.6}
In order to apply the results in section 1, we need to know the 
condition when $F(\r) = \ol\r$ for an irreducible character 
$\r$ of $G^{F^2}$.  
We return to the setting in 2.2, and further assume that 
$F(s) = s\iv$.  Then $F$ acts on $W_s$, preserving 
$W_s^0$ and $\Om_s$.  We denote this 
action by  $\g$, so that $\g^2 = \d$.  
Note that if $\r'$ belongs to $\ol\CM_{s, E}$, then $\ol\r'$
belongs to $\ol\CM_{s\iv, E}$ since $E \in \Irr W^0_s$ is self dual.
Also $F(\r')$ belongs to $\ol\CM_{F(s), F(E)}$.
Hence if $\r$ as above belongs to $\ol\CM_{s, E}$, 
the $\Om_s$-orbit of $E$ turns out to be $F$-stable.  
It follows that $\g$ leaves $\wt\Om_s(E)$
invariant, and induces an action on $\wt\Om_s(E)_{\d}$.
We denote by $\wt\Om_s(E)_{\d}^{\g}$ the set of $\g$-fixed points in  
$\wt\Om_s(E)_{\d}$.
$\g$ acts also on the set $\Om_s^{\d}(E)\wg$.  We denote by 
$\Om_s^{\d}(E)\wg_{-\g}$ the set of $\e \in \Om_s^{\d}(E)\wg$ such that
$\g(\e) = \ol\e$.
We put, for $E \in (\ol{\Irr} W_s^0)^{\g}$,  
\begin{equation*}
\ol\CM^0_{s,E} = \Om_s^{\d}(E)\wg_{-\g}\times \wt\Om_s(E)^{\g}_{\d}.
\end{equation*}
We have the following proposition.
\begin{prop} 
Let $\r_{\e,y}$ be the irreducible character of $G^{F^2}$ corresponding to
$(\e,y)  \in \ol\CM_{s, E}$.
Assume that $F(s) = s\iv$.
\begin{enumerate}
\item   
If the $\Om_s$-orbit of $E$ is not $F$-stable, then 
$F(\r_{\e,y}) \ne \ol\r_{\e,y}$.
\item
Assume that the $\Om_s$-orbit of $E$ is $F$-stable.  
Then $F(\r_{\e,y}) = \ol\r_{\e,y}$ if and only if 
$(\e, y)  \in \ol\CM^0_{s, E}$.
\end{enumerate}
\end{prop}
The proposition will be proved in 2.11 after some preliminaries. 
 First we note that 
\begin{lem}  
For each $N$, we have $F(\vG_N) = \vG_N$, and
$\ol{\vG_N} = \vG_N$.
\end{lem}
\begin{proof} 
The fact that $F(\vG_N) = \vG_N$ can be checked directly
for any $N \in \Fg^F$ since $U^{F^2}$ is $F$-stable and 
$\vL_N$ is also $F$-stable.
On the other hand, it follows from the definition that we have  
$\ol{\vG_N} = \vG_{-N}$.
So, in order to show the lemma, it is enough to see that
$N$ is conjugate to $-N$ under $G^{F^2}$. 
Since $N$ is given by a Jordan normal form, this is reduced 
to the case where $N$ is regular nilpotent.  
Assume that $N$ is a regular nilpotent element given in the 
Jordan normal form.  There exists 
$g = \diag(a, -a, \dots, (-1)^{n-1}a) \in \wt G$ 
such that $gNg\iv = -N$.
Then $g \in G^{F^2}$ if and only if $a \in \Fqq$ and 
$(-1)^ka^n = 1$ with $k = [n/2]$. 
We can set $a = 1$ if $k$ is even, and set $a = -1$ if $k$ is 
odd and $n$ is odd. So, assume that $n$ is even and $k$ is odd, 
i.e., $n = 2k$.  In this case, we may take $a \in \Fqq$ such that 
$a^2 = -1$.  Thus we can always find $g \in G^{F^2}$, and the lemma
follows.
\end{proof}
As a corollary, we have
\begin{cor}  
Let $\r_{\ds_x, E_x} \in \Irr \wt G^{F^2}$ and 
$\r_0 \in \Irr G^{F^2}$ be as in 2.5.  
Assume that
$F(\r_{\ds_x, E_x})|_{G^{F^2}} = \ol{\r_{\ds_x, E_x}}|_{G^{F^2}}$.
Then we have $F(\r_0) = \ol{\r_0}$.
\end{cor}
\begin{proof}
The parametrization of $\Irr G^{F^2}$ in terms of the set
$\ol\CM_{s, E}$ is reduced to the special case where 
$\ol\CM_{s, E}$ is given by (2.5.1)
through the steps (b) and (c) in [S3, 4.5].  
Since the steps (b) and (c)
are compatible with the $F$ action and with taking duals, the 
assertion is reduced to the case of (2.5.1).  In this case,
we have $F(\ol{\vG_N}) = \vG_N$ by Lemma 2.8.  
Note that the $F$-action and taking duals preserve the inner product. 
Since $\r_0$
is the unique irreducible character such that 
\begin{equation*}
\lp \vG_N, \r_0 \rp_{G^{F^2}} = \lp \r_{s_x, E_x},\r_0\rp_{G^{F^2}} = 1,  
\end{equation*}
the corollary follows.
\end{proof}
\begin{lem}  
Assume that the set $\ol\CM_{s, E}$ satisfies the assumption of
Proposition 2.7 (ii).  Take $y \in \wt\Om_s(E)_{\d}$ and assume that 
$(E, y) \lra (x, E_x)$ under the map in
(2.2.2). Then
$F(\r_{\ds_x, E_x})|_{G^{F^2}} = \ol{\r_{\ds_x, E_x}}|_{G^{F^2}}$ 
if and only if $y \in \wt\Om_s(E)^{\g}_{\d}$.
\end{lem}
\begin{proof}
We may choose $\dy \in Z_{G^*}(s)$ as a representative of 
$x \in (\Om_s)_{\d}$, so we may assume that   
$E_x = E$. Then it is known by [S3, (4.5.1)] that 
$\CT_{s_x, E_x} = \CT_{s_y, E}$ corresponds to the set 
$\Om^{\d}_s(E)\wg \times \{y\}$ under the correspondence 
$\CT_{s, E} \lra \ol\CM_{s,E}$.
It is easy to see, for any pair $(s_1, E_1)$, 
 that $F(\r_{\ds_1, E_1}) = \r_{F(\ds_1), F(E_1)}$, where
$F(E_1)$ is the character of $W_{F(s_1)}$ corresponding to $E_1$ 
under the isomorphism $W_{s_1} \simeq W_{F(s_1)}$.
On the other hand, we have $\ol{\r_{\ds_1, E_1}} = \r_{\ds_1\iv, E_1}$ since 
$W_{s_1} = W_{s_1\iv}$ and $E_1$ is self dual.
It follows that $F(\ol{\r_{\ds_y, E}}) = \r_{F(\ds_y\iv), F(E)}$.
By our assumption, $F(s\iv) = s$. Hence we have 
$F(s_y\iv) \in T^*_{\g(y)}$ and 
$F(E) \in (\Irr W^0_s)^{\g(y)\d}$.  This implies that 
$F(\ol{\CT_{s_y, E}}) = \CT_{s_{\g(y)}, F(E)} = 
\CT_{s_{\g(y)}, {}^uE}$ for some $u \in \Om_s$ since 
the $\Om_s$-orbit of $E$ is $F$-stable. 
 Since 
$\CT_{s_{\g(y)},{}^uE} = \CT_{s_y, E}$ if and only if $\g(y) = y$ on
$\wt\Om_s(E)_{\d}$, the lemma is proved.
\end{proof}
\para{2.11.}
We shall prove Proposition 2.7.
The assertion (i) follows from 2.6.  We show (ii).
Take $(\e, y) \in \ol\CM_{s, E}$.  
If $\g(y) \ne y$, then $F(\r_{\e, y}) \ne \ol{\r_{\e, y}}$ by 
Lemma 2.10.
So, assume that $y \in \wt\Om_s(E)_{\d}^{\g}$.
Let $\r_{\ds_x, E_x}$ be the character of $\wt G^{F^2}$ containing 
$\r_{\e,y}$. 
Again by Lemma 2.10, we have $F(\r_{\ds_x, E_x}) = \ol{\r_{\ds_x, E_x}}$.
Let $\r_0 \in G^{F^2}$ be as in 2.5.  Then by Corollary 2.9, 
we have $F(\r_0) = \ol\r_0$.
If we write $\r_{\e,y} = {}^g\r_0$ with $g \in \wt G^{F^2}$, we have
$F(\ol\r_{\e,y}) = {}^{F(g)}\r_0$.
Now the action of $F$ induces an action on $(\wt G^{F^2}/G^{F^2})\wg$ 
which is compatible with the natural pairing 
$\wt G^{F^2}/G^{F^2} \times (\wt G^{F^2}/G^{F^2})\wg \to \Ql^*$.
Then $F$ stabilizes the subgroup $I(\r_{\ds_x, E_x})$.
\par
The arguments in 2.4 shows that the condition 
$F(\r_{\e,y}) = \ol\r_{\e,y}$ is described by investigating 
the action of $F$ on $I(\r_{\ds_x, E_x})$. We follow the notation 
in 2.4.  $I(\r_{\ds_x, E_x})$ is regarded as a subgroup of $\wt Z^{*F^2}$.
If we denote by $\wt\w_{s_x}$ the map $\Om_s^{\d} \to \wt Z^{*F^2}$
obtained as the composite of $\w_{s_x}$ and the inclusion 
$\wt Z^{*F^2}_{s_x} \hra \wt Z^{*F^2}$, then $\wt\w_{s_x}(\Om_s^{\d}(E_x))$
coincides with $I(\r_{\ds_x, E_x})$.  
We note the following.
\par\medskip\noindent
(2.11.1) \ Assume that $x \in (\Om_s)_{\d}$ is $\g$-stable.  Then 
the following diagram commutes. 
\begin{equation*}
\begin{CD}
\Om_s^{\d} @>\wt\w_{s_x}>>  \wt Z^{*F^2}\\
@ V-\g VV                 @VV F = \g V  \\
\Om_s^{\d} @>\wt\w_{s_x}>>  \wt Z^{*F^2},  
\end{CD}
\end{equation*}
where $-\g: \Om_s^{\d} \to \Om_s^{\d}$ is the map defined 
by $z \mapsto \g(z)\iv$.  
\par\medskip
We show (2.11.1).
Take $z \in \Om_s^{\d}$.  Since $F(s) = s\iv$, we have 
\begin{equation*}
\g(\wt\w_{s_x}(z)) = F({}^{g_x}(\ds\iv\dz\ds\dz\iv)) 
       =  {}^{F(g_x)}(\ds F(\dz)\ds\iv F(\dz)\iv).
\end{equation*}
On the other hand, since $x$ is $\g$-stable, $\wt\w_{s_x}$ 
coincides with $\wt\w_{s_{\g(x)}}$, and we have
\begin{equation*}
\wt\w_{s_x}(-\g(z)) = {}^{F(g_x)}(\ds\iv F(\dz)\iv \ds F(\dz))
                    = {}^{F(g_x)}(\ds F(\dz)\ds\iv F(\dz)\iv).
\end{equation*} 
since $\ds\iv F(\dz)\iv \ds F(\dz)$ is in the center $\wt Z^*$
of $\wt G^*$.  Hence (2.11.1) holds.
\par
Now (2.11.1) shows that the $F$-action on $I(\r_{\ds_x, E_x})$ is 
transferred to the $-\g$ action on $\Om_s^{\d}(E)$.
Hence under the parametrization 
$\CT_{s_x, E_x} \lra \Om_s^{\d}(E)\times \{ y\}$ given by 
$\r_{\e,y} \lra (\e, y)$, we see that $F(\r_{\e,y}) = \ol{\r_{\e, y}}$
if and only if $\e$ is $-\g$ stable, i.e., 
$\e \in \Om_s^{\d}(E)\wg_{-\g}$.
This proves the proposition.
\section{Almost characters of $SL_n(\Fqq)$}
\para{3.1.}
We shall parametrize $F^2$-stable irreducible characters of 
$G^{F^{2m}}$ for a sufficiently divisible integer $m$.
Let $s$ be an  $F^2$-stable semisimple element in $G^*$.
We assume that $m$ is large enough so that  
 $F^{2m}$ acts trivially on $W_s^0$ and $\Om_s$.  We denote by
$\ol\CM_{s, E}^{(m)}$ the set parametrizing irreducible 
characters of $G^{F^{2m}}$ corresponding to 
$\ol\CM_{s, E}$ in the previous section.  Hence, 
$\ol\CM_{s, E}^{(m)} = \Om_s(E)\wg\times \Om_s(E)$.
Since the class $s$ is $F^2$-stable, one can define a map
$\d = F^2: W_s \to W_s$ as before.  If the $\Om_s$-orbit 
of $E$ is $F^2$-stable, then $\d$ stabilizes 
$\Om_s(E)$.  
For a pair $(s, E)$ such that $E \in (\ol{\Irr} W_s^0)^{\d}$,
we define a subset $\CM_{s, E}$ of $\ol\CM_{s, E}^{(m)}$ by
\begin{equation*}
\CM_{s, E} = \Om_s(E)_{\d}\wg \times \Om_s(E)^{\d},  
\end{equation*}
where $\Om_s(E)\wg_{\d}$ is the set of $\d$-stable irreducible 
characters in $\Om_s(E)\wg$. 
We denote by $\CE(G^{F^{2m}}, \{s\})^{F^2}$ the subset of 
$F^2$-stable irreducible characters in $\CE(G^{F^{2m}},\{ s\})$.
Then by [S3, (4.6.1)] it is known that under the parametrization 
in (2.2.1) for $G^{F^{2m}}$, we have
\begin{equation*}
\CE(G^{F^{2m}}, \{ s\})^{F^2} = \coprod_{E \in (\ol{\Irr} W_s)^{\d}}
                        \CM_{s, E}.
\end{equation*} 
We denote by $\r_{\e,z}^{(m)}$ the $F^2$-stable irreducible character 
of $G^{F^{2m}}$ corresponding to $(\e,z) \in \CM_{s,E}$.
\par
In the case where $F(s) = s\iv$, 
one can define a map $\g= F: W_{s} \to W_{s}$
preserving $\Om_s$ and $W^0_s$, and such that $\d = \g^2$ as before.  
We denote by $\Om_s(E)^{\g}$ the $\g$-fixed point subgroup of 
$\Om_s(E)$, and by $\Om_s(E)\wg_{-\g}$ the set of 
$\e \in \Om_s(E)\wg$ such that $\g(\e) = \ol\e$.  Then we define a 
subset $\CM_{s,E}^0$ of $\CM_{s, E}$ by 
\begin{equation*}
\CM_{s, E}^0 = \Om_s(E)\wg_{-\g} \times \Om_s(E)^{\g}.
\end{equation*}
\par
The following proposition can be proved in a similar way as 
in Proposition 2.7.
\begin{prop}  
Let $\r^{(m)}_x$ be an $F^2$-stable 
irreducible character of $G^{F^{2m}}$ corresponding to 
$x \in \CM_{s, E}$.  Assume that $F(s) = s\iv$. 
\begin{enumerate}
\item
If the $\Om_s$-orbit of $E$ is not $F$-stable, then 
$F(\r_x^{(m)}) \ne \ol\r^{(m)}_x$.
\item 
Assume that the $\Om_s$-orbit of $E$ is $F$-stable.
Then $F(\r_x^{(m)}) = \ol\r^{(m)}_x$ if and only if 
$x \in \CM_{s, E}^0$.  
\end{enumerate}
\end{prop}
\para{3.3.}
Following [S3, 4.6], we define almost characters of $G^{F^2}$.
For a given $\wt\Om_s(E)_{\d}$, we choose $a = a_E \in \Om_s$ and write
it as $\wt\Om_s(E)_{\d} = \Om_s(E)_{\d}a_E$.
For $x = (\e, z) \in \CM_{ s, E}$ and 
$y = (\e', z'a) \in \ol\CM_{s, E}$, we define  
a pairing $\{ x, y\} \in \Ql^*$ by
\begin{equation*}
\{ x , y\} = |\Om_s(E)^{\d}|\iv \e(z')\e'(z).
\end{equation*}
Then we define a class function $R_x$ of $G^{F^2}$ by 
\begin{equation*}
\tag{3.3.1}
R_x = \sum_{y \in \ol\CM_{s, E}}\{ x,y\}\r_y.
\end{equation*}
$R_x$ are called almost characters of $G^{F^2}$.  Note that 
the definition of the pairing $\{ \ ,\ \}$ depends on the choice
of $a_E \in \wt\Om_s(E)_{\d}$.  If $a_E$ is replaced by $a' = b\iv a_E$
with $b \in \Om_s(E)_{\d}$, then $R_x$ is replaced by 
$\e(b)R_x$.  Hence the almost character $R_x$ is determined uniquely
up to a root of unity multiple.
\par
It is easy to see that (3.3.1) can be converted to the form
\begin{equation*}
\tag{3.3.2}
\r_y = \sum_{x \in \CM_{s, E}}\{ x,y\}\iv R_x.
\end{equation*}
\par
The following result describes the Shintani descent of 
irreducible characters of $G^{F^{2m}}$.  Here we write 
the restriction of $F^2$ on $G^{F^{2m}}$ as $\d$ instead of 
$\s^2$, in connection with the previous section. 
\begin{thm}[{[S3, Theorem 4.7]}]  
Let $\r^{(m)}_x$ be an $F^2$-stable irreducible character of 
$G^{F^{2m}}$ corresponding to $x \in \CM_{s, E}$, and 
choose an extension 
$\wt\r_x^{(m)}$ to $G^{F^{2m}}\lp\d\rp$.  Then 
\begin{equation*}
Sh_{F^{2m}/F^2}(\wt\r^{(m)}_x|_{G^{F^{2m}}\d}) = \mu_xR_x,
\end{equation*}
where $\mu_x$ is a root of unity depending on the extension 
$\wt\r_x^{(m)}$ and on the choice of $a_E$.
\end{thm}
\par
Combining Theorem 3.4 with Proposition 3.2, we have the following 
refinement of Corollary 1.18.
\begin{cor}  
Let $\CM_{s, E}$ be such that $F(s) = s\iv$ and that
the $\Om_s$-orbit of $E$ is $F$-stable.  Then 
\begin{equation*}
 m_2(R_x) = \begin{cases} 
                \z_x    &\quad\text{ if } x \in \CM_{s, E}^0, \\
                0       &\quad\text{ otherwise},   
            \end{cases}
\end{equation*} 
where $\z_x$ is a certain root of unity.
\end{cor} 
The following result describes the action of twisting operators on 
almost characters.  In the special case where $F^2$ acts trivially 
on the center, this was proved by Bonnaf\'e 
[B, Th\'eor\`eme 5.5.4].  We note that this result is also
derived from the property of character sheaves, by making use of
Lusztig's conjecture for $SL_n$, which will be discussed in [S4].
\begin{thm}  
For any $x = (\e,z) \in \CM_{s,E}$, we have
\begin{equation*}
t_1^*(R_x) = \e(z)\iv R_x.
\end{equation*}
\end{thm}
The theorem will be proved in 3.15 after some preliminaries. First we 
recall some general properties of twisting operators.
\begin{lem} 
Let $\vG$ be a connected algebraic group defined over $\Fq$ with 
Frobenius map $F$, and $H$ a connected $F$-stable subgroup 
of $\vG$.  Then the twisting operator
$t_1^*$ commutes with the induction $\Ind_{H^F}^{\vG^F}$.
\end{lem}
\begin{proof}
It is clear that $t_1^*$ commutes with the restriction functor
$\Res_{H^F}^{\vG^F}$.  Moreover $t_1^*$ is an isometry with respect to the
inner product on $H^F$ and $\vG^F$.  The lemma follows from these two facts. 
\end{proof}
\para{3.8}
Let $\vG$ be as in the lemma.  For each integer $m > 0$, we consider 
the group $\vG^{F^m}$, and its semidirect product
$\wt\vG^{F^m} = \vG^{F^m}\lp\s\rp$, where $\s$ is the restriction 
of $F$ on $\vG^{F^m}$,
and $\lp\s\rp$ is the cyclic group of order $m$ with generator $\s$. 
Then the  twisting operator $t_1^*:  C(\vG^F\ssim) \to C(\vG^F\ssim)$
can be lifted to the operator 
\begin{equation*}
\t_1^{*-1} : C(\vG^{F^m}\s\ssim) \to C(\vG^{F^m}\s\ssim)
\end{equation*}
in the following way.
We define a map $\t_1: \vG^{F^m}\s\ssim \to \vG^{F^m}\s\ssim$ by 
$\t_1(x\s) = (x\s)^{1-m}$, and define $\t_1^*$ by its transpose.
It is shown in [S1, Lemma 4.2] that, under the condition 
that $m$ is sufficiently divisible, $\t_1^*$ is an isomorphism and
satisfies the following commutative diagram.
\begin{equation*}
\tag{3.8.1}
\begin{CD}
C(\vG^{F^m}\s\ssim) @>\t_1^*>> C(\vG^{F^m}\s\ssim) \\
@V Sh_{F^m/F} VV               @VV Sh_{F^m/F} V  \\
C(\vG^F\ssim)  @< t_1^* <<     C(\vG^F\ssim).
\end{CD}
\end{equation*}  
We have the following result.
\begin{thm}[{[S1, Theorem 4.7]}]  
Let $\wt\r$ be an extension of an $F$-stable irreducible character 
of $\vG^{F^m}$ to $\wt\vG^{F^m}$.  Then for an appropriate choice 
of (sufficiently divisible) $m$, there exists a root of unity $\la$ 
such that 
\begin{equation*}
\t_1^*(\wt\r|_{\vG^{F^m}\s}) = \la(\wt\r|_{\vG^{F^m}\s}).
\end{equation*}
\end{thm}
\par
The following related result seems to be worth mentioning though
it is not used later.  In [S1], under some condition on $p$, 
the notion of almost characters was established for any connected
algebraic group $\vG$.  Then in view of (3.8.1) together with
Theorem 3.9, we have
\begin{cor}  
For each almost character $R_x$ of $\vG^F$, there exists a root 
of unity $\la_x$ such that 
\begin{equation*}
t_1^*(R_x) = \la_xR_x.
\end{equation*} 
\end{cor}
The following result was proved by Digne and Michel, which holds 
for any connected reductive groups.
\begin{prop}[{[DM]}] 
Let $H$ be an $F$-stable Levi subgroup of a parabolic subgroup 
of a connected reductive group $\vG$.  Then the Lusztig induction 
$R_{H}^{\vG}: C(H^F\ssim) \to C(\vG^F\ssim)$ commutes with 
the twisting operator $t_1^*$.
\end{prop}
\para{3.12.}
We now return to our original setting, and consider $G = SL_n$.
The modified generalized Gelfand-Graev characters were introduced by 
Kawanaka (see [K1]), which is a refinement of generalized Gelfand-Graev 
characters.  The modified generalized Gelfand-Graev characters are 
used in [S3] to parametrize irreducible characters of $SL_n$.  Here
we discuss the action of twisting operators on the modified
generalized Gelfand-Graev characters.  
We follow the notation in 2.5 (but replacing $F^2$ by $F$). 
\par
By [S3, 2.6], we may choose $\Fu$
so that $\Fu$ is $L$-stable. 
Let $A_{\la} = Z_L(\la)/Z^0_L(\la)$.
Then by [S3, 2.7], we have
\begin{equation*}
A_{\la} \simeq A_G(N) =  Z_G(N)/Z^0_G(N).
\end{equation*}
In particular, $A_{\la}$ is an abelian group.
$F$ acts naturally on $A_{\la}$, and we consider 
the quotient group $(A_{\la})_F$ of $A_{\la}$. 
Put 
\begin{equation*}
\ol\CM = (A_{\la})_F \times (A_{\la}^F)\wg.
\end{equation*} 
For each pair $(c, \xi) \in \ol\CM$ one can define a modified 
generalized Gelfand-Graev character $\vG_{c,\xi}$
as follows.  For $c \in A_{\la}$, we choose a representative
$\dc \in Z_L(\la)$.  Then we find $\a_c \in L$ such that
$\a_c\iv F(\a_c) = \dc$.  We define a linear 
map $\la_c: \Fu \to k$ by $\la_c = \la\circ\Ad \a_c\iv$, 
and define a linear 
character $\vL_c = \p\circ\la_c\circ\log$ on $U^F$.
Since $Z_L(\la_c)^F = Z_L(\vL_c)^F$, the linear character
$\vL_c$ can be extended to a linear character on 
$Z_L(\la_c)^FU^F$ trivial on $Z_L^0(\la_c)^F$, which we denote
also by $\vL_c$.
On the other hand, since $A_{\la}$ is abelian, 
we can define a linear character 
$\xi\nat$ of $Z_L(\la_c)^F$ trivial on $Z_L^0(\la_c)^F$ by 
\begin{equation*}
\xi\nat: Z_L(\la_c)^F \to (Z_L(\la_c)/Z_L^0(\la_c))^F \simeq
     A_{\la}^{\dc F} = A_{\la}^F \to \Ql^*,
\end{equation*} 
where the last step is given by $\xi: A_{\la}^F \to \Ql^*$.
We denote by the same symbol $\xi\nat$ the lift of $\xi\nat$ 
to $Z_L(\la_c)^FU^F$ under the homomorphism 
$Z_L(\la_c)^FU^F \to Z_L(\la_c)^F$. Under these setting 
we define $\vG_{c,\xi}$ by 
\begin{equation*}
\vG_{c,\xi} = \Ind_{Z_L(\la_c)^FU^F}^{G^F}
          (\xi\nat\otimes\vL_c).
\end{equation*} 
\para{3.13.}
We choose $m$ large enough so that $F^m$ acts trivially on 
$A_{\la}$.  Replacing $F$ by $F^m$, we have a modified 
generalized Gelfand-Graev character $\vG_{(c,\xi)}^{(m)}$ on 
$G^{F^m}$.    Now the parameter set $\ol\CM$ is replaced by 
$A_{\la} \times (A_{\la})\wg$.
We denote by $\CM$ the subset of 
$A_{\la} \times (A_{\la})\wg$ defined by 
\begin{equation*}
\CM = A_{\la}^F \times (A_{\la})\wg_F,
\end{equation*}
where $(A_{\la})\wg_F$ is the set of $F$-stable irreducible 
characters of $A_{\la}$.
Following [S3, 1.8], we construct, for each $(c, \xi) \in \CM$, 
an $F$-stable modified generalized Gelfand-Graev character 
$\vG^{(m)}_{c,\xi}$, and its extension to 
$G^{F^m}\lp\s\rp$, where $\s = F|_{G^{F^m}}$.   
For $c \in A_{\la}^F$, we choose $\dc \in L^F$. 
We construct the linear character 
$\vL^{(m)}_c$ of $U^{F^m}$ as in 3.12, i.e.,   
we choose $\b_c \in L$
such that $\b_c\iv F^m(\b_c) = \dc$, 
and define $\la_c$
by $\la_c = \la\circ\Ad\b_c\iv$, and put 
$\vL_c^{(m)} = \p_m\circ\la_c\circ\log$, where 
$\p_m = \p\circ\Tr_{{\mathbf F}_{q^m}/\Fq}$.  Put 
$\hat c = \b_cF(\b_c\iv) \in L^{F^m}$. Then  
$\vL_c^{(m)}$ turns out to be $\hat cF$-stable.
\par
On the other hand, it can be checked that $\hat cF$ acts
on $Z_L(\la_c)$ commuting with $F^m$, and that
under the isomorphism 
\begin{equation*}
\ad \b_c\iv : Z_L(\la_c)^{F^m}/Z_L^0(\la_c)^{F^m} 
  \simeq Z_L(\la)^{\dc F^m}/Z_L^0(\la)^{\dc F^m}
  \simeq A_{\la},
\end{equation*}
the action of $\hat cF$ on $Z_L(\la_c)^{F^m}$
is transferred to the action of $F$ on $A_{\la}$. 
Hence if we take $\xi \in (A_{\la})\wg_F$, it produces
an $\hat cF$-stable linear character $\xi\nat$ 
on $Z_L(\la_c)^{F^m}$.
It follows that 
 $\xi\nat\otimes \vL^{(m)}_c$ is $\hat cF$-stable 
for $(c, \xi) \in \CM$, and we conclude that 
$\vG^{(m)}_{c,\xi}$ is $F$-stable.
\par
Put $\hat c_0 = (\hat c\s)^m \in L^{F^m}$.
We note that 
$\hat c_0 \in Z_L(\la_c)^{F^m} = Z_L(\vL_c^{(m)})^{F^m}$.
In fact, since $\vL_c^{(m)}$ is $\hat cF$-stable, 
it is stable by $(\hat c\s)^{m} = \hat c_0$.
We also note that 
\begin{equation*}
\b_c\iv \wh c_0\b_c \equiv \dc\iv \pmod {Z_L^0(\la)^{\dc F^m}}
\end{equation*}
since $(\hat c\s)^m = \b_c F^m(\b_c\iv)$ and 
$\dc = \b_c\iv F^m(\b_c)$.  
In particular, we have 
\begin{equation*}
\tag{3.13.1}
\xi\nat(\hat c_0) = \xi(c\iv).
\end{equation*}
Put $M_c = Z_L(\la_c)^{F^m}$ and $M_c^0 = Z_L^0(\la_c)^{F^m}$.
We consider a subgroup $M_cU^{F^m}\lp\hat c\s\rp$ of 
$G^{F^m}\lp\s\rp$ generated by $M_cU^{F^m}$ and $\hat c\s$. 
Since $\xi\nat \in M_c\wg$ is $\hat cF$-stable, and 
$(\hat c\s)^m = \hat c_0 \in M_c$, $\xi\nat$ may be extended 
to a linear character $\wt\xi\nat$ of $M_c\lp\hat c\s\rp$
in $m$ distinct way.  
The extension $\wt\xi\nat$ is determined by the value 
$\wt\xi\nat(\hat c\s) = \mu_{c,\xi}$, where $\mu_{c,\xi}$ is
any $m$-th root of $\xi\nat(\hat c_0)$.
\par
We fix an extension $\wt\xi\nat$ of $\xi\nat$ to 
$M_c\lp\hat c\s\rp$.
Since $M_cU^{F^m}\lp\hat c\s\rp$ is the semidirect product
of $M_c\lp\hat c\s\rp$ with $U^{F^m}$, $\wt\xi\nat$ may be
regarded as a character of $M_cU^{F^m}\lp\hat c\s\rp$.
On the other hand, since $\wt\vL_c^{(m)}$ is $\hat c\s$-stable,
it can be extended to a linear character on 
$M_cU^{F^m}\lp\hat c\s\rp$ by $\wt\vL_c^{(m)}(\hat c\s) = 1$.
Thus we have a character $\wt\xi\nat\otimes\wt\vL_c^{(m)}$ of
$M_cU^{F^m}\lp\hat c\s\rp$ which is an extension of 
$\xi\nat\otimes\vL_c^{(m)}$ on $M_cU^{F^m}$.
We put
\begin{equation*}
\wt\vG_{c,\xi}^{(m)} = 
   \Ind_{M_cU^{F^m}\lp\hat c\s\rp}^{G^{F^m}\lp\s\rp}
                            (\wt\xi\nat\otimes\wt\vL_c^{(m)}).
\end{equation*}
Then $\wt\vG_{c,\xi}^{(m)}$ gives rise to an extension of $\vG^{(m)}_{c,\xi}$
to $G^{F^m}\lp\s\rp$.  Note that 
$\mu_{c,\xi}\iv \wt\vG^{(m)}_{c,\xi}|_{G^{F^m}\s}$ depends only on
the choice of $(c, \xi)$.
\par
Now we have the following result. 
\begin{prop}  
Let the notations be as above. We have
\begin{equation*} 
\t_1^*(\wt\vG_{c,\xi}^{(m)}|_{G^{F^m}\s}) 
         = \xi(c)(\wt\vG_{c,\xi}^{(m)}|_{G^{F^m}\s})
\end{equation*}
for an appropriate choice of (sufficiently divisible) $m$.
\end{prop}
\begin{proof}
The following proof is an analogy of the argument in 
[S1, Corollary 5.10].
Put $H = LU$.  Then $H$ is an $F$-stable  connected subgroup of 
$G$. 
For each $(c,\xi) \in  \CM$, we denote by $\th$ the linear character
$\xi\nat\otimes\vL_c^{(m)}$ of $M_cU^{F^m}$, and 
by $\wt\th$ its extension $\wt\xi\nat\otimes \wt\vL_c^{(m)}$ 
of $M_cU^{F^m}\lp\wh c\s\rp$.  We put $\wt H^{F^m} = H^{F^m}\lp\s\rp$,
$V = M_cU^{F^m}$, and $\wt V = M_cU^{F^m}\lp\hat c\s\rp$.
We consider the induced characters
\begin{equation*}
\r_{c,\xi} = 
           \Ind_{V}^{H^{F^m}}\th,
\quad 
\wt\r_{c,\xi} = \Ind_{\wt V}^{\wt H^{F^m}}\wt\th.
\end{equation*}
Then $\r_{c,\xi}$ is an $F$-stable character of $H^{F^m}$, and 
$\wt\r_{c,\xi}$ is an extension of $\r_{c,\xi}$ to $\wt H^{F^m}$.
Moreover, $\r_{c,\xi}$ is irreducible by [S3, Lemma 1.7].
Note that 
\begin{equation*}
\wt\vG^{(m)}_{c,\xi}|_{G^{F^m}\s} = 
        \Ind_{H^{F^m}\s}^{G^{F^m}\s}(\wt\r_{c,\xi}|_{H^{F^m}\s}). 
\end{equation*}
In order to prove the proposition, we have only to show the following 
formula since $\t_1^*$ commutes with the induction
$\Ind_{H^{F^m}\s}^{G^{F^m}\s}$ by Lemma 3.7 and (3.8.1).
\begin{equation*}
\tag{3.14.1}
\t_1^*(\wt\r_{c,\xi}|_{H^{F^m}\s}) =
                \xi(c)(\wt\r_{c,\xi}|_{H^{F^m}\s}).
\end{equation*} 
We show (3.14.1).  We choose $m$ so that $m$ is a multiple of some 
fixed integer $A$, where $A$ is divisible by $|A_{\la}|p$, 
and that $m-1$ is prime to the order of 
$\wt H^{F^m}$.  The existence of such $m$ is 
shown in [S1, Lemma 4.8].  Then the map 
$f: \wt H^{F^m} \to \wt H^{F^m}, g \mapsto g^{1-m}$ 
is a bijection, and $\t_1$ is obtained by restricting $f$ to 
$H^{F^m}\s$.  Since 
$f$ stabilizes the conjugacy classes, it induces an isomorphism
$f^*: C(\wt H^{F^m}\ssim) \to C(\wt H^{F^m}\ssim)$.  The map 
$f^*$ stabilizes the space $C(\wt V\ssim)$, and we denote by 
$f^*_{\wt V}$ the restriction of $f^*$ on  $\wt V$. 
We note that 
\par\medskip\noindent
(3.14.2) \ $f_{\wt V}^*(\wt\th)$ is a linear character of 
$\wt V$ such  that  $f^*_{\wt V}(\wt\th)|_V = \th$.
\par\medskip
In fact,  since $f$ induces a homomorphism on $\wt H^{F^m}$ modulo
the commutator subgroup, 
$f_{\wt V}^*$ maps linear characters to linear characters.
We show that the restriction of $f_{\wt V}^*(\wt\th)$ on $V$
coincides with $\th$.  
Since $\la\circ\log : U \to k$ is a homomorphism of algebraic groups
and $m$ is divisble by $p$, $\la(g^m) = 0$ for $g \in U$.  This
implies 
that $\vL_c^{(m)}(g^{1-m}) = \vL_c^{(m)}(g)$ for $g \in U^{F^m}$.
On the other hand, since $m$ is divisible by $|A_{\la}|$, 
$\xi\nat(g^{1-m}) = \xi\nat(g)$ for $g \in M_c$.  
It follows that $\th(g^{1-m}) = \th(g)$ for any 
$g \in V = M_cU^{F^m}$, and the claim follows.
\par
Now it is easy to see that $f^*$ commutes with the induction
\begin{equation*}
\Ind_{\wt V}^{\wt H^{F^m}}: C(\wt V\ssim) \to C(\wt H^{F^m}\ssim).
\end{equation*}
Thus $f^*(\wt\r_{c,\xi})$ is also an extension of $\r_{c,\xi}$ to
$\wt H^{F^m}$.
Now the extensions of $\th$ to $\wt\th$ is characterized by the
value $\wt\th(\hat c\s)$, and it determines the extension 
$\wt\r_{c,\xi}$.
Since $f(\hat c\s) = (\hat c\s)^{1-m} = \hat c\s\cdot \hat c_0\iv$, 
we see that 
\begin{align*}
f^*_{\wt V}(\wt\xi\nat\otimes \wt\vL_c^{(m)})(\hat c\s)
     &= \xi\nat(\hat c_0\iv)\cdot 
           \wt\xi\nat\otimes \wt\vL_c^{(m)}(\hat c\s) \\
     &= \xi(c)\cdot \wt\xi\nat\otimes \wt\vL_c^{(m)}(\hat c\s)
\end{align*}
by (3.13.1). 
This proves (3.14.1), and so the proposition follows.
\end{proof}  
\para{3.15.}
We are in a position to prove Theorem 3.6.  
We apply the previous results to our situation by replacing $F$ by 
$F^2$.  By [S3, 4.5], the parametrization of 
irreducible characters, $\r_{\e',z'} \lra (\e', z') \in \ol\CM_{s,E}$
are divided into three steps.  Accordingly, the parametrization of
almost characters $R_{\e,z} \lra (\e,z) \in \CM_{s,E}$ are divided 
similarly.  The cases (b) and (c) in [loc. cit.] are reduced to the
case (a) via Harish-Chadra induction and Lusztig induction.  Since 
the twisting operator commutes with Lusztig induction by 
Proposition 3.11, the proof of Theorem 3.6 is reduced to the case (a),
i.e., the case where $(s, E)$ satisfies the condition (2.5.1).
\par
So, assume that $(s, E)$ is as above.  In this case, irreducible 
characters belonging to $\ol\CM_{s, E}$ and almost characters 
belonging to $\CM_{s, E}$ are characterized by modified generalized 
Gelfand-Graev characters as follows. 
Let $\r_{\ds, E}$ be an irreducible character of $\wt G^{F^2}$ for 
some $\ds \in \wt G^*$ such that $\pi(\ds) = s$, and $N$ be a
nilpotent element such that the nilpotent orbit $\CO_N$ containing $N$
coincides with the orbit $\CO$ associated to $\r_{\ds, E}$ (see 
e.g., [S3, 2.9])
Let $A_{\la}$ be the finite group given in 3.12. For a certain
quotient group $\ol A_{\la}$ of $A_{\la}$ with $F^2$-action, 
we put 
\begin{equation*}
\CM_{s, N} = \ol A_{\la}^{\d} \times (\ol A_{\la})\wg_{\d},
\end{equation*} 
where $\d$ is the action of $F^2$ on $\ol A_{\la}$ and on 
$(\ol A_{\la})\wg$ as before.
Let $(\CT_{s, E}^{(m)})^{F^2}$ be the set of $F^2$-stable 
irreducible characters of $G^{F^{2m}}$ belonging to $\CM_{s, E}$. 
Then there exists a parametrization 
$\CM_{s, N} \lra (\CT_{s, E}^{(m)})^{F^2}$ via 
$(c,\xi) \lra \r_{c,\xi}^{(m)}$ satisfying the following 
properties.
Put $\wt\CM_{s, N} = A^{\d}_{\la} \times (\ol A_{\la})\wg_{\d}$.
Since $\wt\CM_{s, N}$ is a subset of 
$A_{\la}^{\d} \times (A_{\la})\wg_{\d}$, one can define an
$F^2$-stable character $\vG^{(m)}_{c,\xi}$ of $G^{F^{2m}}$ 
for each pair $(c, \xi) \in \CM_{s, N}$ (see 3.13).
Let $\vf: \wt\CM_{s, N} \to \CM_{s, N}$ be the natural projection.
Then for $(c, \xi) \in \wt\CM_{s, N}$ and $(c',\xi') \in \CM_{s, N}$, 
we have the following (cf. [S3, Corollary 2.21]).
\begin{equation*}
\tag{3.15.1}
\lp\vG^{(m)}_{c,\xi}, \r^{(m)}_{c',\xi'}\rp = 
     \begin{cases}
          1 &\quad\text{ if } \vf(c, \xi) = (c', \xi'), \\
          0 &\quad\text{ otherwise.}
     \end{cases} 
\end{equation*}
Assume that $\vf(c,\xi) = (c', \xi')$, i.e., $\xi = \xi'$ and 
$c'$ is the image of $c$ under the map $A_{\la} \to \ol A_{\la}$.
Let $\wt\vG_{c,\xi}^{(m)}$ and $\wt\r^{(m)}_{c',\xi'}$ be 
extensions of $\vG^{(m)}_{c,\xi}$ and $\r^{(m)}_{c',\xi'}$ to
$G^{F^{2m}}\lp\d\rp$, respectively. 
Now by Proposition 3.14, $\wt\vG_{c,\xi}^{(m)}|_{G^{F^{2m}}\d}$ 
is an eigenfunction  
for $\t^*_1$ with eigenvalue $\xi(c)$.  Since $\r_{c',\xi'}^{(m)}$
occurs in the decomposition of $\vG_{c,\xi}^{(m)}$ with 
multiplicity 1 by (3.15.1), by applying Theorem 3.9 we see that 
\begin{equation*}
\tag{3.15.2}
\t_1^*(\wt\r_{c',\xi'}^{(m)}|_{G^{F^{2m}}\d}) =
         \xi'(c')(\wt\r_{c',\xi'}^{(m)}|_{G^{F^{2m}}\d}).
\end{equation*} 
(Note that $\xi(c) = \xi'(c')$).
The set $(\CT_{s,E}^{(m)})^{F^2}$ is also parametrized by 
the set $\CM_{s, E}$ via $\r_{\e,z}^{(m)} \lra (\e,z) \in \CM_{s,E}$.
They are related to each other through the bijection 
$\CM_{s,N} \lra \CM_{s, E}$, $(c',\xi') \lra (\e, z)$ which satisfies 
the condition that $\xi'(c') = \e(z)$.
Thus we have 
\begin{equation*}
\tag{3.15.3}
\t_1^*(\wt\r_{\e,z}^{(m)}|_{G^{F^{2m}}\d}) = 
        \e(z)(\wt\r_{\e,z}^{(m)}|_{G^{F^{2m}}\d}).
\end{equation*}
Now the theorem follows from Theorem 3.4, in view of 
the commutativity  of $t_1^{*-1}$ and $\t_1^*$ given in (3.8.1).
This completes the proof of the theorem.
\par\medskip
\section{ Determination of $m_2(\r_{\ds, E}|_{G^{F^2}})$} 
\para{4.1.}
Assume that $s$ is $F^2$-stable, and 
the pair $(s, T^*)$ is given as in Section 2.
Let $\r_{\ds, E}$ be an irreducible character of $\wt G^{F^2}$
as in 2.2.
In this section, we shall compute the value 
$m_2(\r_{\ds, E}|_{G^{F^2}})$.
Now $\r_{\ds, E}$ is given as  
\begin{equation*}
\tag{4.1.1}
\r_{\ds, E} = \ve_{\wt G^*}\ve_{Z_{\wt G}(\ds)}|W_{\ds}|\iv
         \sum_{w \in W_{\ds}}\Tr(w\d, \wt E)R_{\wt T^*_{w}}(\ds),
\end{equation*}
where $\ve_H = (-1)^{\Fqq-\text{rank}(H)}$ for any reductive group
$H$, and $\wt E$ is a certain  
extension of $E \in (\ol\Irr W_{\ds})^{\d}$ to $W_{\ds}\lp\d\rp$.
Therefore we compute the value $m_2(R_{\wt T^*_w}(\ds))$ 
for each $w \in W_s$.
\par
We consider the isomorphism 
$\wt Z^{*F^2} \simeq (\wt G^{F^2}/G^{F^2})\wg$ as in (2.4.1),
and a similar one by replacing $F^2$ by $F$.  
By the property of the dual torus, 
we have the following commutative diagram.
\begin{equation*}
\begin{CD}
\wt Z^{*F^2} @> \sim >> (\wt G^{F^2}/G^{F^2})\wg \\
@V N_{F^2/F}VV              @VV \Res V    \\
\wt Z^{*F}   @> \sim >>  (\wt G^F/G^F)\wg, 
\end{CD}
\tag{4.1.2}
\end{equation*}
where $\Res$ is the restriction of the character of 
$\wt G^{F^2}/G^{F^2}$ on $\wt G^F/G^F$, and 
$N_{F^2/F}$ is the norm map $z \to zF(z)$.  
The norm map is also described as in 1.1. By using this,
it is easy to see that $\Ker N_{F^2/F}$ coincides with 
the subset $\{ z\iv F(z) \mid z \in \wt Z^{*F^2}\}$, and so
$\wt Z^{*F}$ can be identified with $(\wt Z^{*F^2})_F$ 
via the map $zF(z) \lra z$ for $z \in (\wt Z^{*F^2})_F$.
\par
First we note the following general fact.
\begin{lem} 
Let $\x$ be a class function of $\wt G^{F^2}$.  Then 
\begin{equation*}
m_2(\x|_{G^{F^2}}) = \sum_{\th \in (\wt G^{F}/G^{F})\wg}
                       m_2(\x\otimes \wt\th),
\end{equation*}
where $\wt\th$ is a character of $\wt G^{F^2}/G^{F^2}$,
 regarded as a linear character of $\wt G^{F^2}$, 
which is an extension of $\th$ via the inclusion
$\wt G^F/G^F \hra \wt G^{F^2}/G^{F^2}$. 
\end{lem}
\begin{proof}  By the Frobenius reciprocity, we have
\begin{align*}
\lp\x|_{G^{F^2}}, \Ind_{G^F}^{G^{F^2}}1\rp &=
                    \lp\x, \Ind_{G^F}^{\wt G^{F^2}}1\rp \\
    &= \lp\x, \Ind_{\wt G^F}^{\wt G^{F^2}}
                \bigl(\Ind_{G^F}^{\wt G^F}1\bigr)\rp \\
    &= \lp\x,   \sum_{\th \in (\wt G^{F}/G^{F})\wg} 
                 \Ind_{\wt G^F}^{\wt G^{F^2}}\th\rp.
\end{align*}    
But for any linear character $\wt\th$ of $\wt G^{F^2}$ such that
$\wt\th|_{\wt G^{F}} = \th$, we have            
\begin{equation*}
\lp\x, \Ind_{\wt G^F}^{\wt G^{F^2}}\th\rp = 
\lp \x\otimes\wt\th\iv, \Ind_{\wt G^F}^{\wt G^{F^2}}1\rp 
   = m_2(\x\otimes\wt\th\iv).
\end{equation*}
Thus the lemma is proved.
\end{proof}
By applying the above formula to the class function 
$R_{\wt T^*_w}(\ds)$ of $\wt G^{F^2}$, 
\begin{lem} 
We have
\begin{equation*}
\tag{4.3.1}
m_2(R_{\wt T^*_{w}}(\ds)|_{G^{F^2}}) 
       = \sum_{z \in (\wt Z^{*F^2})_F}
                     m_2(R_{\wt T^*_w}(\ds\dz)),
\end{equation*}
where $\dz$ is a representative of $z$ in $\wt Z^{*F^2}$.
\end{lem}
\begin{proof}
By 4.1, $(\wt Z^{*F^2})_F$ is isomorphic to
$(\wt G^F/G^F)\wg$.  We denote by $\th$ the character of $\wt G^F/G^F$
corresponding to $z \in (\wt Z^{*F^2})_F$.
Then by (4.1.2), the representative $\dz \in \wt Z^{*F^2}$ 
corresponds to a linear
character $\wt\th$ of $\wt G^{F^2}/G^{F^2}$, which is an 
extension of $\th$.  Now it is known that 
$R_{\wt T^*_w}(\ds)\otimes\wt\th = R_{\wt T^*_w}(\ds\dz)$.
Hence the lemma follows from Lemma~4.2.
\end{proof}
We have the following proposition.
\begin{prop}  
Let $s$ be an element in $T_w^*$ such that $F^2(s) = s$.
\begin{enumerate}
\item
Assume that the $G^{*F^2}$-orbit of $s$ does not contain $s'$ such that 
$F(s') = {s'}\iv$.   
Then $m_2(R_{\wt T^*_w}(\ds)|_{G^{F^2}}) = 0$ 
for any $\ds \in \wt T_w^*$ such that $\pi(\ds) = s$.
\item
Assume that $F(s) = s\iv$.  Then
there exists $\ds \in \wt T_w^*$ such that $\pi(\ds) = s$ and that 
$F(\ds) = \ds\iv$, and we have 
\begin{equation*}
\tag{4.4.1}
m_2(R_{\wt T^*_w}(\ds)|_{G^{F^2}}) = 
       \sum_{x \in \Om_s^{-\g}}m_2(R_{\wt T^*_w}(\ds z_x)),
\end{equation*}
where $\Om_s^{-\g} = \{ x \in \Om_s \mid \g(x) = x\iv\}$ is the 
subgroup of $\Om_s^{\d}$, and  
$z_x \in \wt Z^{*F^2}$ is a representative of an element in
$(\wt Z^{*F^2})_F$ such that $z_xF(z_x) = \w_s(x)$ under the map 
$\w_s: \Om_s^{\d} \to \wt Z_{\ds}^{*F^2} \subset \wt Z^{*F^2}$ (see 2.4). 
\end{enumerate}
\end{prop}
\begin{proof}
First we show (i).  
It is known that $R_{\wt T^*_w}(\ds)|_{G^{F^2}}$ coincides with 
the Deligne-Lusztig character $R_{T^*_w}(s)$ of $G^{F^2}$.
Then by [L3, 2.7 (a)], we have $m_2(R_{T^*_w}(s)) = 0$.  
(Note that in [loc. cit.], it is assumed that 
the center of $G$ is connected.  However, the above fact 
holds without this assumption, by 2.3 and 2.6 (b) in [loc. cit.]).  
Thus (i) holds.  
\par
Next we show (ii).  Assume that $F(s) = s\iv$.  
Take $\ds_1 \in \wt T^*$ such that $\pi(\ds_1) = s$.
Then there exists some $z \in  \wt Z^*$ such that $F(\ds_1) = \ds_1\iv z$ for 
some $z \in \wt Z^*$.  Take $z_1 \in \wt Z^* \simeq \BG_m$ such that 
$z = z_1F(z_1) = z_1^{q+1}$, and put $\ds = \ds_1z_1$.  
Then $\pi(\ds) = s$ and $F(\ds) = \ds\iv$ as asserted.
\par
Take $\ds$ as above, and consider the formula (4.3.1).  Again by 
[L3, Lemma 2.8], we may only consider, in the sum of the right hand side
of (4.3.1),   
 $z \in (\wt Z^{*F^2})_F$ such that $F(\ds\dz)$ is conjugate to 
$(\ds\dz)\iv$ in $\wt G^*$.  Here we note that 
\par\medskip\noindent
(4.4.2) \ $F(\ds\dz)$ is conjugate to $(\ds\dz)\iv$ 
if and only if there exists $x \in \Om_s^{-\g}$ such that 
$\dz F(\dz) = \w_s(x)$.
\par\medskip
We show (4.4.2).  Assume that $x \in \Om_s^{-\g}$, and let 
$\dx$ be an element in $\wt G^*$ such that $\pi(\dx)$ is a
representative of $x$.  Then by (2.11.1), 
$\w_s(x) = \ds\iv\dx \ds\dx\iv \in \wt Z^{*F^2}$ is $\g$-stable, i.e.,
$\w_s(x) \in \wt Z^{*F}$.  Hence by (4.1.2) 
there exists $\dz \in \wt Z^{*F^2}$ such that
$\w_s(x) = \dz F(\dz)$. It follows that 
$\ds\dz = \dx \ds\dx\iv F(\dz)\iv$, and we have 
$F(\ds\dz) = \dx\iv (\ds\dz)\iv\dx$.  This shows that $F(\ds\dz)$ 
is conjugate to $(\ds\dz)\iv$ in $\wt G^*$.
Conversely, assume that $F(\ds\dz)$ is conjugate to $(\ds\dz)\iv$ in
$\wt G^*$.  Then there exists $\dx \in \wt G^*$ such that 
$F(\ds\dz) = \dx\iv (\ds\dz)\iv \dx$.  Clearly 
$\pi(\dx) \in Z_{G^*}(s)$, and its image in $\Om_s$ determines 
an element $x \in \Om_s$.  Since $\ds(\dz F(\dz)) = \dx \ds\dx\iv$,
we have $\dz F(\dz) \in \wt Z^{*F^2}_s$.  Moreover, $\dz F(\dz)$ is
$\g$-stable.  Hence by (2.4.2) and (2.11.1), we see that
$x \in \Om_s^{-\g}$.  This proves (4.4.2).
\par
Since $\w_s(x) \in \wt Z^{*F}$, $\dz \in \wt Z^{*F^2}$ such that
$\dz F(\dz) = \w_s(x)$ has 
a unique image on $(\wt Z^{*F^2})_F \simeq \wt Z^{*F}$ by 
(4.1.2).  We choose $z_x$ from such $\dz$ for each $x$.  Then 
the formula (4.4.1) is immediate from (4.4.2). 
\end{proof}
By using Lusztig's formula in [L3], we shall compute the right hand 
side of (4.4.1) explicitly.  We show 
\begin{lem}  
Under the notation in Proposition 4.4 (ii), we have
\begin{equation*}
m_2(R_{\wt T^*_w}(\ds z_x)) = \sharp\{ u \in W_{\ds} \mid w = u({}^{x\g}u)\}.
\end{equation*}
\end{lem}
\begin{proof}
Take $\dx \in N_{\wt G^*}(\wt T^*)$ whose image in $W$ is a representative
of $x \in \Om_s^{-\g}$.   
We note that one can choose $\dx$ such that $F(\dx) = \dx\iv$.  In
fact, take any $x' \in N_{\wt G^*}(\wt T^*)$ in the inverse image 
of $x$.  Since $\g(x) = x\iv$,
we have $x' F(x') = t \in \wt T^*$. We can find $t_1 \in \wt T^*$ 
such that $t_1\iv F^2(t_1) = t$.  Then $\dx = t_1x'F(t_1)\iv$ satisfies 
the required condition. 
\par
Take $g \in \wt G^*$ such that $g\iv F(g) =  \dx$.
Put $s' = {}^{g}(\ds z_x)$, $\wt T' = {}^g\wt T^*$ 
and $W' = N_{\wt G^*}(\wt T')/\wt T'$.
Then $F(s') = {s'}\iv$, $F(\wt T') = \wt T'$, 
and $s' \in \wt T'$.
Moreover, we have $g\iv F^2(g) = \dx F(\dx) = 1$ and so
$g \in \wt G^{*F^2}$. 
We have an isomorphism 
$f: W \isom W'$ via $\ad g$, and we see that the pair
$(s', \wt T'_{f(w)})$ is $\wt G^{*F^2}$-conjugate to the pair
$(\ds z_x, \wt T^*_w)$, where $\wt T'_{f(w)}$ is an $F^2$-stable maximal
torus obtained from $\wt T'$ by twisting by $f(w) \in W'$.  
\par
It follows that 
\begin{equation*}
\tag{4.5.1}
R_{\wt T^*_{w}}(\ds z_x) = R_{\wt T'_{f(w)}}(s').
\end{equation*}
If we put $W'_{s'} = \{ w' \in W' \mid w'(s') = s' \}$,  
$f$ induces an isomorphism $W_{\ds} \isom W'_{s'}$.  Now
It is known by [L3, Lemma 2.8, (b)] that 
\begin{equation*}
\tag{4.5.2}
m_2(R_{\wt T'_{f(w)}}(s')) = \sharp\{ y \in W'_{s'} \mid f(w) = yF(y)\}. 
\end{equation*}  
Hence, by (4.5.1) and (4.5.2), we have 
\begin{align*}
m_2(R_{\wt T_w}(\ds z_x)) 
            &= \sharp\{ u \in W_{\ds} \mid f(w) = f(u)F(f(u))\}  \\
            &= \sharp\{ u \in W_{\ds} \mid w = u({}^{x\g}u) \}, 
\end{align*}
since $F\circ f = f\circ \dx F$.  This proves the lemma.  
\end{proof}
We are in a position to determine $m_2(\r_{\ds, E}|_{G^{F^2}})$. 
\begin{thm}  
Let $\r_{\ds, E}$ be an irreducible character of $\wt G^{F^2}$
as before, and put $s = \pi(\ds)$. 
\begin{enumerate}
\item 
If $s$ is not $G^{*F^2}$-conjugate to $s'$ such that 
$F(s') = {s'}\iv$, then 
$m_2(\r_{\ds, E}|_{G^{F^2}}) = 0$. 
\item
Assume that $F(s) = s\iv$. 
Then we have
\begin{equation*}
m_2(\r_{\ds,E}|_{G^{F^2}}) = \begin{cases}
                 |\Om_s^{-\g}(E)| 
   &\quad\text{ if there exists $x \in \Om_s^{-\g}$ such that 
                ${}^{x\g}E = E$ }, \\
                 0 &\quad\text{ otherwise},
                           \end{cases}
\end{equation*}
where $\Om_s^{-\g}(E)$ is the stabilizer of $E$ in $\Om_s^{-\g}$. 
\end{enumerate} 
\end{thm}
\begin{proof} 
$\r_{\ds, E}$ is given as in (4.1.1).  Thus 
(i) is immediate from Proposition 4.4 (i).  We show (ii).  So,
assume that $F(s) = s\iv$.
Since $R_{\wt T^*_w(\ds)}|_{G^{F^2}}$ does not depend on the
choice of a representative $\ds$ of $s$, we may assume that
$\ds$ satisfies the property that $F(\ds) = \ds\iv$.
Then $\ve_{\wt G*}\ve_{Z_{\wt G^*}(s)} = 1$ by [L3, 1.5 (b)].
Hence by (4.4.1) together with Lemma 4.5, we have 
\begin{align*}
m_2(\r_{\ds, E}|_{G^{F^2}}) 
     &= |W_{\ds}|\iv\sum_{w \in W_{\ds}}
        Tr(w\d, \wt E) m_2(R_{\wt T_w^*}(\ds)|_{G^{F^2}}) \\
     &= |W_{\ds}|\iv\sum_{w \in W_{\ds}}\Tr(w\d, \wt E)
           \sum_{x \in \Om_s^{-\g}}m_2(R_{\wt T_w^*}(\ds z_x)) \\
     & = |W_{\ds}|\iv\sum_{x \in \Om_s^{-\g}}\sum_{w \in W_{\ds}}
     \Tr(w\d, \wt E)\sharp\{ u \in W_{\ds} \mid w = u({}^{x\g}u)\}.
\end{align*}
\par
Now by Lemma 2.11 in [L3] (see also the formula in the proof of 
Proposition 2.13 there), one can write 
\begin{equation*}
\sum_{E' \in (W_{\ds})\wg_{x\g}}\Tr(w\d, \wt E') = 
        \sharp\{ u \in W_{\ds} \mid w = u({}^{x\g}u)\}, 
\end{equation*}
where $(W_{\ds})\wg_{x\g}$ is the set of $x\g$-stable characters of 
$W_{\ds}$, and the extension $\wt E'$ of $E'$ is chosen to be
realized over $\BQ$.   
(Note that $(x\g)^2 = \d$ since $x \in \Om_s^{-\g}$).
It follows that 
\begin{equation*}
m_2(\r_{\ds,E}|_{G^{F^2}}) = 
       \sum_{x \in \Om_s^{-\g}}\sum_{E' \in (W_{\ds})\wg_{x\g}}
      |W_{\ds}|\iv\sum_{w \in W_{\ds}}\Tr(w\d, \wt E)\Tr(w\d, \wt E').    
\end{equation*}
But since 
\begin{equation*}
|W_{\ds}|\iv\sum_{w \in W_{\ds}}\Tr(w\d, \wt E)\Tr(w\d, \wt E') =
          \begin{cases}
                   1  &\quad\text{ if } \wt E = \wt E', \\
                   0  &\quad \text{ if }  E \ne E', 
          \end{cases}
\end{equation*}
(here the extension $\wt E$ is chosen to be over $\BQ$, see [L1, 3.2]),
we have
\begin{equation*}
m_2(\r_{\ds,E}|_{G^{F^2}}) = \sharp\{ x \in \Om_s^{-\g} \mid 
                      \text{ ${}^{x\g}E = E$ } \}. 
\end{equation*}
If there exists $x_1$ such that ${}^{x_1\g}E = E$, then 
$\{ x \in \Om_s^{-\g} \mid {}^{x\g}E = E\} = \Om_s^{-\g}(E)x_1$.
Thus the theorem is proved.
\end{proof}
\par
We shall apply the formula in the theorem to the case 
$\r_{\ds_x, E'}$.    First we note that 
\begin{lem}  
Assume that $F(s) = s\iv$. 
Let $s_y$ be an element corresponding to $ y \in (\Om_s)_{\d}$.
Then the $G^{*F^2}$-class  of $s_y$ contains 
an element $s'$ such that 
$F(s') = {s'}\iv$ if and only if there exists $u \in \Om_s$ such that
$uF(u)$ gives a representative of $y$ in $\Om_s$.
\end{lem}
\begin{proof}
Let $\dy \in Z_{G^*}(s)$ be a representative of $y$.  
Let $s'$ be an element contained in the $G^{*F^2}$-class of 
$s_y$.  Then $s'$ can be obtained as $s' = {}^gs$ for some  
$g \in G^*$ such that $g\iv F^2(g) = \dy$.
It is easy to see that 
$F(s') = {s'}\iv$ if and only if $g\iv F(g) \in Z_{G^*}(s)$. 
Hence $y = uF(u)$ in $\Om_s$ if we put $u$ the image of 
$g\iv F(g)$ in $\Om_s$.
\end{proof}
\para{4.8.}
We prepare a notation. Let $s$ be a semisimple element 
such that $F(s) = s\iv$, and $E \in \Irr W^0_s$ such that
the $\Om_s$-orbit of $E$ is $F$-stable. 
We define a subset $\wt\Om_s(E)_{\d}^+$ 
(resp. $\Om_s(E)_{\d}^+$) of $\wt\Om_s(E)_{\d}$ 
(resp. of $\Om_s(E)_{\d}$) by 
\begin{align*}
\wt\Om_s(E)^+_{\d} &= \text{ the image of }
  \{ u\g(u) \mid u \in \Om_s, {}^{u\g}E = E \} 
  \text{ into } \wt\Om_s(E)_{\d}, \\
\Om_s(E)^+_{\d} &= \text{ the image of }
   \{ v\g(v) \mid v \in \Om_s(E)\} \text{ into } \Om_s(E)_{\d}.
\end{align*}
Then we can see that there exists $a_E \in \Om_s$ such that
\begin{equation*}
\tag{4.8.1}
\wt\Om_s(E)^+_{\d} = \Om_s(E)^+_{\d}a_E.
\end{equation*}
In fact, since the $\Om_s$-orbit of $E$ is $\g$-stable, there
exists $b \in \Om_s$ such that ${}^{b\g}E = E$.
Then $a_E = b\g(b)$ is contained in $\wt \Om_s(E)$, and we have
$\wt\Om_s(E) = \Om_s(E)a_E$.
(4.8.1) follows from this.
\par
As a corollary to Theorem 4.6, we have the following.
\begin{cor} 
Assume that $s$ is semisimple in $G^*$ such that $F(s) =
s\iv$, and that $E \in (\ol\Irr W_s^0)^{\d}$.
 Take $y \in \wt\Om_s(E)_{\d}$ and 
let $(E, y) \lra (x, E_x)$ be as in (2.2.2).
Then we have
\begin{enumerate}
\item
If the $\Om_s$-orbit of $E$ is not $F$-stable, then 
$m_2(\r_{\ds_x, E_x}) = 0$.
\item  
Assume that the $\Om_s$-orbit of $E$ is $F$-stable.  Then 
\begin{equation*}
m_2(\r_{\ds_x, E_x}) = \begin{cases}
            |\Om_s^{-\g}(E)| 
        &\quad\text{ if } y \in \wt\Om_s(E)^+_{\d},  \\
            0 &\quad\text{ otherwise}.
                      \end{cases} 
\end{equation*}
\end{enumerate}
\end{cor}
\begin{proof}
If  $m_2(\r_{\ds_x, E_x}) \ne 0$,  
then $s_x$ is $G^{F^2}$-conjugate to some $s'$ such that 
$F(s') = {s'}\iv$.  Since $x$ and $y$ are in the same class
in $(\Om_s)_{\d}$, there exists $u \in \Om_s$ such that 
$uF(u) = y$ by Lemma 4.7.  Let $\du \in Z_{G^*}(s)$
be a representative of $u$.  Then there exists $g \in G^*$
such that $g\iv F(g) = \du$.   We see that $g\iv F^2(g) = \du F(\du)$
is a representative of $y$.  Hence we may assume 
$s_x = s_y = {}^gs$ and $E_x = E$.  
Then $\ad g\iv$ gives rise to an isomorphism 
$W^0_{s_x} \isom W^0_s$, and $\ad g\iv$ sends $F, F^2$ to $uF, yF^2$,
respectively.  Moreover, $E'' \in \Irr W_{s_x}^0$ in (2.2.2) is mapped
to $E' = E \in \Irr W_s^0$.   
Then by Theorem 4.6, $m_2(\r_{s_x, E_x}) \ne 0$ is equivalent 
to the condition that there exists $h \in \Om^{-\g}_s(E)$ such 
that ${}^{hu\g}E = E$.  
In particular, the $\Om_s$-orbit of $E$ is $F$-stable. 
Since $hu\g(hu) = u\g(u) = y$, this is
equivalent to $y \in \wt\Om_s(E)^+_{\d}$.
This proves the corollary.
\end{proof}
\par\medskip
\section{Determination of $m_2(\r)$ for $\r \in \Irr SL_n(\Fqq)$}
\para{5.1.}
In this section, we shall determine $m_2(\r)$ for all 
irreducible characters of $G^{F^2}$.  Our strategy is to
compute $m_2$ for almost characters of $G^{F^2}$ first, 
and then derive the formula for $m_2(\r)$ from it. 
\par
First we prepare some notation.  
Let $s$ be an $F^2$-stable semisimple element in $G^*$, and 
$E$ an $F^2$-stable irreducible character of $W^0_s$. We recall
two sets 
$\ol\CM_{s,E} = \Om_s^{\d}(E)\wg \times \wt\Om_s(E)_{\d}$ and 
$\CM_{s, E} = \Om_s(E)\wg_{\d} \times \Om_s(E)^{\d}$.
Assuming that $F(s) = s\iv$ and that the $\Om_s$-orbit of 
$E$ is $F$-stable,  we define subsets 
$
\Om_s^{\d}(E)\wg_- \subset \Om_s^{\d}(E)\wg_{-\g} 
                    \subset \Om_s^{\d}(E)\wg
$
by 
\begin{align*}
\Om_s^{\d}(E)\wg_{-\g} &= 
      \{ \th \in \Om_s^{\d}(E)\wg \mid \g(\th) = \th\iv\},  \\ 
\Om_s^{\d}(E)\wg_-     &= 
      \{ \th\iv\g(\th) \mid \th \in \Om_s^{\d}(E)\wg\}.
\end{align*}
We also consider subsets 
$
\wt\Om_s(E)^+_{\d}  \subset \wt\Om_s(E)_{\d}^{\g} \subset 
              \wt\Om_s(E)_{\d},
$
where
\begin{equation*}
\wt\Om_s(E)_{\d}^{\g} = \{ u \in \wt\Om_s(E)_{\d} \mid \g(u) = u\},
\end{equation*}
and $\wt\Om_s(E)_{\d}^+$ is defined as in 4.8.
We define subsets 
$\Om_s(E)^+_{\d} \subset \Om_s(E)_{\d}^{\g} \subset \Om_s(E)_{\d}$ 
in a similar way as above.
\par
Put 
\begin{equation*}
|\Om_s(E)^{\d}| = t, \quad |\Om_s(E)^{\g}| = d, \quad
|\Om_s(E)^{-\g}| = d'.
\end{equation*}
Then we see easily that
\begin{align*}
|\Om_s^{\d}(E)\wg| &= |\Om_s(E)_{\d}| = t, \\
|\Om_s^{\d}(E)\wg_{\g}| &= |\Om_s(E)_{\d}^{\g}| = d, \\
|\Om_s^{\d}(E)\wg_{-\g}| &= |\Om_s(E)_{\d}^{-\g}| = d'.
\end{align*}
Since $\Om_s(E)^{\d}$ is a cyclic group, $\Om_s(E)^{\d}$ is
written as a product of $\Om_s(E)^{\g}$ and $\Om_s(E)^{-\g}$.
If $t = |\Om_s(E)^{\d}|$ is even, $\Om_s(E)^{\d}$ contains a
unique element of order 2.  In that case  
$\Om_s(E)^{\g}$ and $\Om_s(E)^{-\g}$ has a non-trivial 
intersection, and so $t = dd'/2$.  If $t$ is odd, then 
$\Om_s(E)^{\d} = \Om_s(E)^{\g}\times \Om_s(E)^{-\g}$, and so
$t = dd'$.
\par
There is a surjective homomorphism 
$\Om_s^{\d}(E)\wg \to \Om_s^{\d}(E)\wg_-$ given by 
$\th \to \th\iv F(\th)$, whose kernel is given by
$\Om_s^{\d}(E)\wg_{\g}$.  It follows that 
$\Om_s^{\d}(E)\wg_-$ is a subgroup of $\Om_s^{\d}(E)\wg_{-\g}$
of order $t/d$.  
Hence we have
\begin{equation*} 
\tag{5.1.1}
[\Om_s^{\d}(E)\wg_{-\g} : \Om_s^{\d}(E)\wg_-] 
          = \begin{cases}
                 1 &\quad\text{ if } t = dd', \\
                 2 &\quad\text{ if } t = dd'/2.
             \end{cases}
\end{equation*}
\par
Similarly, we have a surjective homomorphism
$\Om_s(E)_{\d} \to \Om_s(E)^+_{\d}$ given by $z \mapsto zF(z)$ with 
kernel $\Om_s(E)^{-\g}_{\d}$.  It follows that $\Om_s(E)_{\d}^+$ is a
subgroup of $\Om_s(E)_{\d}^{\g}$ of degree $t/d'$.  Hence we have
\begin{equation*}
\tag{5.1.2}
[\Om_s(E)^{\g}_{\d}:\Om_s(E)^+_{\d}] = \begin{cases}
                      1 &\quad\text{ if } t = dd', \\
                      2 &\quad\text{ if } t = dd'/2.
                                  \end{cases}
\end{equation*}
\par
The following result describes the values of $m_2$ for almost 
characters of $G^{F^2}$.
\begin{thm} 
Assume that $s$ is a semisimple element in $G^{*F^2}$, and 
$E$ is an irreducible character of $W_s^0$ such that $\Om_s$-orbit
of $E$ is $F^2$-stable.  Let $R_{\e, z}$ be an almost character 
associated to $(\e,z) \in \CM_{s, E}$.  Then
\begin{enumerate}
\item
Assume that $s$ is not $G^*$-conjugate to an element $s'$ such that 
$F(s') = {s'}\iv$.  Then $m_2(R_{\e, z}) = 0$.
\item
Assume that $F(s) = s\iv$.  If the $\Om_s$-orbit of $E$ is
not $F$-stable, then $m_2(R_{\e,z}) = 0$.
\item
Assume that $F(s) = s\iv$, and that the $\Om_s$-orbit of $E$ is
$F$-stable. 
\begin{enumerate}
\item
Assume that $|\Om_s(E)^{\d}|$ is odd.  Then we have
\begin{equation*}
m_2(R_{\e,z}) = \begin{cases}
                   1 &\quad\text{ if $\e \in \Om_s(E)\wg_{-\g}$
                                    and $z \in \Om_s(E)^{\g}$}, \\
                   0 &\quad\text{ otherwise.}
                 \end{cases} 
\end{equation*}
\item
Assume that $|\Om_s(E)^{\d}|$ is even.  Then we have 
\begin{equation*}
m_2(R_{\e,z}) = \begin{cases}
        1 &\quad\text{ if $\e \in \Om_s(E)\wg_{-\g}$  
           and $z \in \Om_s(E)^+$}, \\ 
        \ve &\quad\text{ if $\e \in \Om_s(E)\wg_{-\g}$  
           and $z \in \Om_s(E)^{\g} - \Om_s(E)^+$}, \\
        0   &\quad\text{ otherwise, }
               \end{cases}                 
\end{equation*}
where $\ve = c_2(\r_{1, z''_0}) = \pm 1$ for any 
$z_0''\in \wt\Om_s(E)^{\g}_{\d} - \wt\Om_s(E)^+_{\d}$. 
\end{enumerate}
\end{enumerate}
\end{thm}
\begin{proof}
We show (i).  By Theorem 4.6 (i), 
$m_2(\r_{\ds_x, E_x}|_{G^{F^2}}) = 0$ for any $x \in (\Om_s)_{\d}$.
It follows that $m_2(\r_{\e',z'}) = 0$ for any 
$(\e',z') \in \ol\CM_{s, E}$.  Hence $m_2(R_{\e,z}) = 0$ for any
$(\e,z) \in \CM_{s,E}$.
\par
A similar proof works for the assertion (ii) 
since $m_2(\r_{\ds_x, E_x}) = 0$ for 
any $x \in (\Om_s)_{\d}$ by Corollary 4.9 (i).  
\par
We show (iii) by computing $m_2(R_{\e, z})$ for $(\e, z) \in \CM_{s, E}$.
By Theorem 3.6, we have 
$m_2(R_{\e,z}) = \e(z)\iv m_2(t_1^{*-1} R_{\e,z})$.  By  
definition of $R_{\e,z}$, together with Corollary 1.11, 
applied to the case where $r = 2$, we have
\begin{equation*}
\tag{5.2.1}
m_2(R_{\e,z}) = \e(z)\iv \sum_{(\e',z'a) \in \ol\CM_{s,E}}
                  \{ (\e,z),(\e', z')\}c_2(\r_{\e',z'a}).
\end{equation*}
with $a = a_E$.
Moreover by Corollary 1.16 (applied to the case where $m = 1$, 
see also [K2, Theorem 2.1.3]), 
together with Proposition 2.7, we have
\begin{equation*}
\tag{5.2.2}
c_2(\r_{\e',z'a}) = \begin{cases}
                     \pm 1 &\quad\text{ if } 
               (\e',z') \in \Om_s^{\d}(E)\wg_{-\g} \times \Om_s(E)^{\g}_{\d}, \\
                      0    &\quad\text{ otherwise.} 
                    \end{cases}
\end{equation*}
We note the following.
\par\medskip\noindent
(5.2.3) \ Assume that $z' \in \Om_s(E)_{\d}^+$. 
  Then $c_2(\r_{\e',z'a}) = 1$ for any 
 $\e' \in \Om_s^{\d}(E)\wg_{-\g}$.
\par\medskip
We show (5.2.3).  Let $(x, E_x)$ corresponding to $(E, z')$ via $f$
in (2.2.2).  Then we have 
$m_2(\r_{\ds_x, E_x}|_{G^{F^2}}) = |\Om_s(E)^{-\g}|$ by Corollary 4.9. 
Since the twisting operator $t_1^*$ acts trivially on $\r_{\ds_x,E_x}$, 
we have $c_2(\r_{\ds_x. E_x}|_{G^{F^2}}) = |\Om_s(E)^{-\g}|$  
by Corollary 1.11.  On the other hand, 
$\r_{\ds_x, E_x}$ can be decomposed as in (2.2.3).  Hence by (5.2.2),  
we have
\begin{equation*}
c_2(\r_{\ds_x,E_x}|_{G^{F^2}}) = \sum_{\e' \in \Om_s^{\d}(E)\wg_{-\g}}
                c_2(\r_{\e',z'a}).
\end{equation*}
Since $|\Om_s(E)^{-\g}| = |\Om_s^{\d}(E)\wg_{-\g}| = d'$, we can 
conclude that $c_2(\r_{\e',z'a}) = 1$, and (5.2.3) follows.
\par
We now compute $m_2(R_{\e,z})$.  
In view of (5.2.2), the formula (5.2.1) can be written as 
\begin{equation*}
\tag{5.2.4}
m_2(R_{\e,z}) = \e(z)\iv|\Om_s^{\d}(E)|\iv
     \sum_{(\e',z') \in \Om_s^{\d}(E)\wg_{-\g} \times \Om_s(E)^{\g}_{\d}}
            \e(z')\e'(z)c_2(\r_{\e',z'a}).
\end{equation*} 
First consider the case where
$t$ is odd, i.e., the case where $t = dd'$.
Then by (5.1.2), we have 
$\Om_s(E)^{\g}_{\d} = \Om_s(E)^+_{\d}$.  It follows by (5.2.3) that 
\begin{equation*}
m_2(R_{\e,z}) = \e(z)\iv|\Om_s^{\d}(E)|\iv
  \sum_{(\e',z') \in \Om_s^{\d}(E)\wg_{-\g} \times \Om_s(E)^{\g}_{\d}}
                     \e(z')\e'(z).
\end{equation*} 
This implies that $m_2(R_{\e,z}) = 0$ unless $\e$ is trivial on
$\Om_s(E)_{\d}^{\g}$ and $z \in \Om_s^{\d}(E)$ is such that
$\e'(z) = 1$ for any $\e' \in \Om_s^{\d}(E)\wg_{-\g}$.
But since $\Om_s^{\d}(E)\wg_{-\g} = \Om_s^{\d}(E)\wg_-$, the condition 
for $z$ is equivalent to the condition that $z \in \Om_s(E)^{\g}$.
Similarly, since $\Om_s(E)^{\g}_{\d} = \Om_s(E)_{\d}^+$, the condition 
for $\e$ is equivalent to the condition that 
$\e \in \Om_s(E)\wg_{-\g}$.
Now assume that $\e \in \Om_s^{\d}(E)\wg_{-\g}$ and 
$z \in \Om_s(E)^{\g}$.
Since $|\Om_s^{\d}(E)| = |\Om_s^{\d}(E)\wg_{-\g}| 
         \times |\Om_s(E)^{\g}_{\d}|$, and $\e(z) = 1$, 
(5.2.4) implies 
that $m_2(R_{\e,z}) = 1$.  This proves (a) of (iii).
\par
Next we consider the case where $t$ is even, i.e., the case where 
$t = dd'/2$.  In this case, $\Om_s^{\d}(E)\wg_-$ is an index two
subgroup of $\Om_s^{\d}(E)\wg_{-\g}$, and $\Om_s(E)^+_{\d}$ is an 
index two subgroup of $\Om_s(E)^{\g}_{\d}$.  We fix 
$\e'_0 \in \Om_s^{\d}(E)\wg_{-\g} - \Om_s^{\d}(E)\wg_-$ and
$z_0' \in \Om_s(E)^{\g}_{\d} - \Om_s(E)^+_{\d}$.
Then by using (5.2.3), (5.2.4) can be written as 
\begin{equation*}
\tag{5.2.5}
m_2(R_{\e,z}) = \e(z)\iv|\Om_s^{\d}(E)|\iv
\sum_{(\e',z') \in \Om_s^{\d}(E)\wg_- \times \Om_s(E)^+_{\d}}
         \e(z')\e'(z)A_{\e',z'}, 
\end{equation*}
where 
\begin{equation*}
A_{\e',z'} = 
1 + \e_0'(z) +  \e(z_0')c_2(\r_{\e',z'z_0'a}) + 
         \e(z_0')\e'_0(z)c_2(\r_{\e'\e_0', z'z_0'a}).
\end{equation*}
It is known by Corollary 3.5 that 
\begin{equation*}
\tag{5.2.6}
|m_2(R_{\e,z})| = \begin{cases}
                     1 &\quad\text{ if } (\e,z) \in 
                          \Om_s(E)\wg_{-\g} \times \Om_s(E)^{\g}, \\
                     0 &\quad\text{ otherwise.}
                   \end{cases}  
\end{equation*}
\par
We now assume that $m_2(R_{\e,z}) \ne 0$. 
Hence $\e \in \Om_s(E)\wg_{-\g}$ and $z \in \Om_s(E)^{\g}$.
We note that $\e \in \Om_s(E)\wg$ is contained in 
$\Om_s(E)\wg_{-\g}$ if and only if 
$\e$ is trivial on $\Om_s(E)^+_{\d}$.  
Similarly, $\e' \in \Om_s^{\d}(E)\wg$ is contained in 
$\Om_s^{\d}(E)\wg_-$ if and only if $\e'$ is trivial on 
$\Om_s(E)^{\g}$. 
In particular,
we have $\e(z') = \e'(z) = 1$ for any 
$(\e',z') \in \Om_s^{\d}(E)\wg_- \times \Om_s(E)^+_{\d}$. 
in the sum in (5.2.5).
Since $\e(z), \e'_0(z), \e(z_0')$ take values $\pm 1$, 
we see that $m_2(R_{\e,z}) \in \BQ$.  This implies 
that $m_2(R_{\e,z}) = \pm 1$ by (5.2.6).
\par
We shall consider the two cases, whether $z \in \Om_s(E)^{\g}$ 
is contained in $\in \Om_s(E)^+$ or not. 
First assume that $z \in \Om_s(E)^+$.
Then $\e(z) = 1$, $\e_0'(z) = 1$ and $\e(z_0') = \pm 1$.
Since
$|\Om_s^{\d}(E)\wg_-|\times |\Om_s(E)^+_{\d}| = t/2$, 
it follows from (5.2.5) that
\begin{equation*}
m_2(R_{\e,z})t = t + \sum_{(\e',z')}\e(z_0')(c_2(\r_{\e',z'z_0'a}) +
             c_2(\r_{\e'e_0', z'z_0'a})).
\end{equation*}
Let $C$ be the sum part of this formula.  Then we have 
$-t \le C \le t$.  Since $m_2(R_{\e,z}) = \pm 1$, this forces that
$C = 0$, and we have $m_2(R_{\e,z}) = 1$.
\par
Next assume that $z \in \Om_s(E)^{\g} - \Om_s(E)^+$.  Since
$z\iv z_0' \in \Om_s(E)^+_{\d}$, we have 
$\e(z\iv z'_0) = 1$. Moreover $\e_0'(z) = -1$.
Hence by (5.2.5), we can write 
\begin{equation*}
m_2(R_{\e,z})t = \sum_{(\e',z')}
            (c_2(\r_{\e',z'z_0'a}) - c_2(\r_{\e'\e_0', z'z_0'a})). 
\end{equation*}
But since
\begin{equation*}
\sum_{(\e',z')}|c_2(\r_{\e',z'z_0'a}) - c_2(\r_{\e'\e_0', z'z_0'a})|
   \le t = |m_2(R_{\e,z})t|,
\end{equation*}
we see that $c_2(\r_{\e',z'z_0'a}) = -c_2(\r_{\e'\e_0', z'z_0'a})$ 
has a common value
for any $(\e',z')$, which coincides with 
$c_2(\r_{1, z_0'a}) = -c_2(\r_{\e_0',z'_0a})$.  This implies 
that $m_2(R_{\e,z}) = c_2(\r_{1, z_0'a}) = \ve$. 
By putting  $z_0'' = z_0'a$, we obtain the theorem. 
\end{proof}
We can now easily translate Theorem 5.2 to the form 
on $m_2(\r)$ for irreducible characters $\r$.
\begin{thm} 
Assume that $s$ is a semisimple element in $G^{*F^2}$, and 
that $E \in\Irr W_s^0$ is such that the $\Om_s$-orbit of $E$ 
is $F^2$-stable.  Let $\r_{\e',\z''}$ be an irreducible character
of $G^{F^2}$ associated to $(\e',z'') \in \ol\CM_{s,E}$.  
Then
\begin{enumerate}
\item
Assume that $s$ is not $G^*$-conjugate to an element $s'$ 
such that $F(s') = {s'}\iv$.  Then $m_2(\r_{\e',z''}) = 0$. 
\item
Assume that $F(s) = s\iv$.  If the $\Om_s$-orbit of $E$ is not
$F$-stable, then $m_2(\r_{\e',z''}) = 0$.
\item
Assume that $F(s) = s\iv$ and that the $\Om_s$-orbit of 
$E$ is $F$-stable. 
\begin{enumerate}
\item
Assume that $|\Om_s(E)^{\d}|$ is odd.  Then we have
\begin{equation*}
m_2(\r_{\e',z''}) = \begin{cases}
                       1  &\quad\text{ if  
        $\e' \in \Om_s^{\d}(E)\wg_{-\g}$ and $z'' \in
                       \wt\Om_s(E)_{\d}^{\g}$,}     \\ 
                       0   &\quad\text{ otherwise. }
                    \end{cases}
\end{equation*}
\item
Assume that $|\Om_s(E)^{\d}|$ is even.  Then we have
\begin{equation*}
m_2(\r_{\e',z''}) = \begin{cases}
                     1 + \ve &\quad\text{ if  
       $\e' \in \Om_s^{\d}(E)\wg_-$ and $z'' \in \wt\Om_s(E)_{\d}^+$} \\
                     1 - \ve &\quad\text{ if 
       $\e' \in \Om_s^{\d}(E)\wg_{-\g} - \Om_s^{\d}(E)\wg_-$ and
              $z'' \in \wt\Om_s(E)^+_{\d}$} \\
                     0       &\quad\text{ otherwise,}
                    \end{cases}
\end{equation*} 
where $\ve = c_2(\r_{1, z_0''}) = \pm 1$ for any 
$z_0'' \in \wt\Om_s(E)^{\g}_{\d} -\wt\Om_s(E)^+_{\d}$.
\end{enumerate}  
\end{enumerate}
\end{thm}
\begin{proof}
The assertion (i) and (ii) are already shown in the proof of 
Theorem 5.2.  We show (iii).  First assume that $|\Om_s(E)^{\d}|$
is odd. Then (3.3.2) implies, in view of Theorem 5.2, that
\begin{equation*}
m_2(\r_{\e',z'a}) = |\Om_s(E)^{\d}|\iv
        \sum_{ (\e,z) \in \Om_s(E)\wg_{-\g} \times \Om_s(E)^{\g}}
               \e(z')\iv\e'(z)\iv.
\end{equation*}  
It follows that $m_2(\r_{\e',z'a}) = 0$ unless $\e'$ is trivial on 
$\Om_s(E)^{\g}$, and $\e(z') = 1$ for any $\e \in \Om_s(E)\wg_{-\g}$, 
and in which  case $m_2(\r_{\e',z'a}) = 1$.
But this condition is equivalent to the condition that 
$\e' \in \Om_s^{\d}(E)\wg_{-\g}$ and $z' \in \Om_s(E)^+_{\d}$.  
By replacing $z'a$ by $z''$, we obtain (a). 
\par
Next assume that $|\Om_s(E)^{\d}|$ is even.
Let us fix $z_0 \in \Om_s(E)^{\g} - \Om_s(E)^+$.  Again by 
Theorem 5.2, we have
\begin{equation*}
m_2(\r_{\e',z'a}) = |\Om_s(E)^{\d}|\iv
       \sum_{(\e,z) \in \Om_s(E)\wg_{-\g} \times \Om_s(E)^+}
               \e(z')\iv\e'(z)\iv(1 + \e'(z_0)\iv\ve).
\end{equation*} 
It follows that $m_2(\r_{\e',z'a}) = 0$ unless $\e'$ is trivial 
on $\Om_s^{\d}(E)^{+}$ and $\e(z') = 1$ for any $\e \in \Om_s(E)\wg_{-\g}$,
and in which case $m_2(\r_{\e',z'a}) = 1 + \e'(z_0)\iv\ve$.
The condition for $z'$ is the same as before, and 
$\e'$ is trivial on $\Om_s^{\d}(E)^{+}$ if and only if 
$\e' \in \Om_s^{\d}(E)\wg_{-\g}$.  Moreover, 
\begin{equation*}
\e'(z_0) = \begin{cases}
           1 &\quad\text{ if } \e' \in \Om_s^{\d}(E)\wg_-, \\
          -1 &\quad\text{ if } \e' \in \Om_s^{\d}(E)\wg_{-\g} 
                                        - \Om_s^{\d}(E)\wg_-.
           \end{cases}
\end{equation*}
Hence (b) holds, and the theorem is proved. 
\end{proof}
\remark{5.4.}  
In [L3], Lusztig gave a uniform description of $m_2(\r)$ for
any irreducible character $\r$ of $G^{F^2}$ in the case where 
$G$ is a connected reductive group with connected center.  He 
expects that his formulation will be extended also to the 
disconnected center case.  
We shall compare our results with Lusztig's conjectural description.
Take $(\e,z) \in \ol\CM_{s, E}$.  
For $z \in \wt\Om_s(E)_{\d}$, take a representative 
$\dz \in \Om_s(E)$ of $z$, and put 
\begin{equation*}
\sqrt z = \text{the image of }
  \{ y \in \Om_s \mid {}^{y\g}E = E, (y\g)^2 = \dz\d\} 
  \text{ into } \Om_s(E)_{\d}.
\end{equation*}
Then $\Om_s(E)^{\d}$ acts on $\sqrt z$ by the $F$-twisted conjugation.
We denote by $\sqrt z$ the corresponding permutation representation
also. Let $[\e: \sqrt z]$ the multiplicity of $\e$ in this permutation
representation.    
Now assume that $F(s) = s\iv$ and that the
$\Om_s$-orbit of $E$ is $F$-stable.  Then we have the following.
\par\medskip\noindent
(5.4.1) \
Assume that 
$(\e,z) \in \Om_s^{\d}(E)\wg_- \times \wt\Om_s(E)^+_{\d}$.  Then we have
\begin{equation*}
[\e: \sqrt z] = \begin{cases}
          1  &\quad\text{ if  $|\Om_s(E)^{\d}|$  is odd},  \\
          2  &\quad\text{ if $|\Om_s(E)^{\d}|$ is even}.  \\
              \end{cases}         
\end{equation*}
If $(\e,z) \notin \Om_s^{\d}(E)\wg_- \times \wt\Om_s(E)^+_{\d}$, 
then we have $[\e : \sqrt z] = 0$.
\par\medskip
In fact, in our setting, $\wt\Om_s(E)^+_{\d}$ is the set of 
$z \in \wt\Om_s(E)_{\d}$ such that $\sqrt z \ne \emptyset$.  
Hence if $z \notin \wt\Om_s(E)^+_{\d}$, then $\sqrt z = \emptyset$, 
and so $[\e:\sqrt z] = 0$.  If $z \in \wt\Om_s(E)^+_{\d}$, then
$\sqrt z = \Om_s(E)_{\d}^{-\g}y$ for some $y \in \sqrt z$.
Let $\x$ be the character of the representation $\sqrt z$. Then
$\x(u) = |\Om_s(E)_{\d}^{-\g}|$ if $u \in \Om_s(E)^{\g}$ and 
$\x(u) = 0$ otherwise.  It follows that
\begin{equation*}
[\e: \sqrt z] = \begin{cases}
         |\Om_s(E)_{\d}^{-\g}||\Om_s(E)^{\g}|/|\Om_s(E)^{\d}| 
             &\quad\text{ if } \e \in \Om_s^{\d}(E)\wg_-, \\
         0   &\quad\text{ otherwise.}
                 \end{cases}    
\end{equation*}
(5.4.1) follows from this.
\par
In the connected center case, $[\e: \sqrt z]$ gives the value
$m_2(\r)$.  In our situation, by comparing with Theorem 5.3, we see
that $[\e:\sqrt z]$ coincides with $m_2(\r_{\e,z})$ if 
$c_2(\r_{1,z_0}) = 1$ for 
$z_0 \in \wt\Om_s(E)^{\g}_{\d} - \wt\Om_s(E)^+_{\d}$. 
\par
\vspace{1cm}


\begin{thebibliography}{[DJ]}
\par
\bibitem [B] {B} C. Bonnaf\'e; Op\'erateurs de torsion dans 
$SL_n(\FF_q)$ et dans $SU_n(\FF_q)$, Bull. Soc. Math. France
{\bf 128} (2000), 309 - 345.  
\par
\bibitem [DM] {DM} F. Digne et J. Michel; Descente de Shintani 
des caracteres de Deligne-Lusztig, CRAS {\bf 291} (1980), 651-653.
\par
\bibitem [H] {H} A. Henderson; Spherical functions of the symmetric
space $G(\FF_{q^2})/G(\FF_q)$, Represent. Theory, {\bf 5} (2001), 
581 - 614.
\par
\bibitem [K1] {K1} N. Kawanaka; Shintani lifting and generalized 
Gelfand-Graev representations, in ``The Arcata conference on 
Representations of of Finite Groups'' Proceedings of Symposia
in Pure Math., Vol. {\bf 47}, pp. 147 - 163, Amer. Math. Soc., 
Providence, RI, 1987.
\par
\bibitem [K2] {K2} N. Kawanaka; 
On subfield symmetric spaces over a finite filed, Osaka J. Math. 
{\bf 28} (1991) 759 - 791.
\par
\bibitem [KM] {KM} N. Kawanaka and H. Matsuyama; A twisted version of 
Frobenius-Schur indicator and multiplicity free permutation
representations, Hokkaido Math. J. {\bf 19} (1990), 495 - 508. 
\par
\bibitem[L1]{L1}G. Lusztig;  \lq\lq Characters of reductive groups
  over a
  Finite field,'' Ann. of Math. Studies, Vol.107, Princeton
  Univ. Press,   Princeton, 1984.
\par
\bibitem [L2] {L} G. Lusztig; On the representations of reductive 
groups with disconnected centre, Asterisque {\bf 168} (1988), 157 - 168.
\par
\bibitem [L3] {L} G. Lusztig;
$G(F_q)$-invariants in irreducible $G(F_{q^2})$-modules, 
Represent. Theory, {\bf 4} (2000), 446 - 465.
\par
\bibitem [S1] {S1} T. Shoji; Shintani descent for algebraic groups 
over a finite field, I, J. Alg., {\bf 145} (1992), 468 - 524.

\bibitem [S2] {S2}T. Shoji; Character sheaves and almost characters of
reductive groups, Adv. in Math. {\bf 111} (1995), 244 - 313, II, Adv. 
in Math. {\bf 111} (1995), 314 - 354. 
\par
\bibitem [S3] {S3} T. Shoji; Shintani descent for special linear
  groups, J. of Algebra {\bf 199} (1998), 175 - 228. 
\par
\bibitem [S4] {S4} T. Shoji; Lusztig's conjecture for finite
special linear groups, in preparation.
\end{thebibliography}
\end{document}